\input amstex
\input amsppt.sty

\define\C{\Bbb C}
\define\Z{\Bbb Z}
\define\Mat{\operatorname{Mat}(m,\C)}
\define\tht{\thetag}
\define\diag{\operatorname{diag}}
\define\M{\Cal M}

\define\wt{\widetilde}
\define\A{\Cal A}
\redefine\B{\Cal B}
\define\Y{\Cal Y}
\define\s{\Cal S}
\define\ze{\zeta}
\define\R{\Cal R}

\hoffset 1cm
\NoBlackBoxes
\TagsOnRight
\NoRunningHeads

\topmatter
\title Isomonodromy transformations of linear systems of difference equations
\endtitle
\author Alexei Borodin
\endauthor
\abstract
We introduce and study ``isomonodromy'' transformations of the matrix linear difference equation $Y(z+1)=A(z)Y(z)$ with polynomial (or rational) $A(z)$. Our main result is a construction of an isomonodromy action of $\Z^{m(n+1)-1}$ on the space of coefficients $A(z)$ (here $m$ is the size of matrices and $n$ is the degree of $A(z)$). The (birational) action of certain rank $n$ subgroups can be described by difference analogs of the classical Schlesinger equations, and we prove that for generic initial conditions these difference Schlesinger equations have a unique solution. We also show that both the classical Schlesinger equations and the Schlesinger transformations known in the isomonodromy theory, can be obtained as limits of our action in two different limit regimes.

Similarly to the continuous case, for $m=n=2$ the difference Schlesinger equations and their $q$-analogs yield discrete Painlev\'e equations; examples include dPII, dPIV, dPV, and $q$-PVI.
\endabstract
\endtopmatter
\head
Introduction
\endhead

In recent years there has been considerable interest in analyzing a certain class of discrete probabilistic models which in appropriate limits converge to well-known models of Random Matrix Theory. The sources of these models are quite diverse, they include Combinatorics, Representation Theory, Percolation Theory, Random Growth Processes, tiling models and others. 

One quantity of interest in both discrete models and their random matrix limits is the {\it gap probability} -- the probability of having no particles in a given set. It is known, due to  works of many people, see \cite{JMMS}, \cite{Me}, \cite{TW}, \cite{P}, \cite{HI}, \cite{BD}, that in the continuous (random matrix type) setup these probabilities can be expressed through a solution of an associated isomonodromy problem for a linear system of differential equations with rational coefficients. 

The goal of this paper is to develop a general theory of ``isomonodromy'' transformations for linear systems of {\it difference} equations with rational coefficients. This subject is of interest in its own right. As an application of the theory, we show in a subsequent publication that the gap probabilities in the discrete models mentioned above are expressible through solutions of isomonodromy problems for such systems of difference equations. In the case of one-interval gap probability this has been done (in a different language) in \cite{Bor}, \cite{BB}. One example of the probabilistic models in question can be found at the end of this Introduction.

Consider a matrix linear difference equation 
$$
Y(z+1)=A(z)Y(z).
\tag 1
$$
Here 
$$
A(z)=A_0z^n+A_1z^{n-1}+\dots+A_n, \qquad A_i\in\Mat,
$$ 
is a matrix polynomial and $Y:\C\to\Mat$ is a matrix meromorphic function. 
\footnote{Changing $Y(z)$ to $(\Gamma(z))^k Y(z)$  readily reduces
a rational $A(z)$ to a polynomial one.} 
We assume that the eigenvalues of $A_0$ are nonzero and that their ratios are not real.
Then, without loss of generality, we may assume that $A_0$ is diagonal. 

It is a fundamental result proved by Birkhoff in 1911, that the equation \tht{1} has two canonical meromorphic solutions $Y^l(z)$ and $Y^r(z)$, which are holomorphic and invertible for $\Re z\ll 0 $ and $\Re z\gg 0$ respectively, and whose asymptotics at $z=\infty$ in any left (right) half-plane has a certain form. Birkhoff further showed that the ratio 
$$
P(z)=(Y^r(z))^{-1}Y^l(z),
$$
which must be periodic for obvious reasons, is, in fact, a rational function in $\exp(2\pi i z)$. This rational function has just as many constants involved as there are matrix elements in $A_1,\dots,A_n$. Let us call $P(z)$ the {\it monodromy matrix} of $\tht{1}$.

 Other results of Birkhoff show that for any periodic matrix $P$ of a specific form, there exists an equation of the form \tht{1} with prescribed $A_0$, which has $P$ as the monodromy matrix. Furthermore, if two equations with coefficients $A(z)$ and $\wt A(z)$, $\wt A_0=A_0$, have the same monodromy matrix, then there exists a rational matrix $R(z)$ such that
$$
\wt A(z)=R(z+1)A(z)R^{-1}(z).
\tag 2
$$

The first result of this paper is a construction, for generic $A(z)$, of a homomorphism of $\Z^{m(n+1)-1}$ into the group of invertible rational matrix functions, such that
the transformation \tht{2} for any $R(z)$ in the image, does not change the monodromy matrix.

If we denote by $a_1,\dots,a_{mn}$ the roots of the equation $\det A(z)=0$ (called {\it eigenvalues of $A(z)$}) and by $d_1,\dots,d_n$ certain uniquely defined exponents of the asymptotic behavior of a canonical solution $Y(z)$ of \tht{1} at $z=\infty$, then
the action of $\Z^{m(n+1)-1}$ is uniquely defined by integral shifts of $\{a_{i}\}$ and $\{d_j\}$ with the total sum of all shifts equal to zero. (We assume that $a_i-a_j\notin\Z$ and $d_i-d_j\notin\Z$ for any $i\ne j$.)

The matrices $R(z)$ depend rationally on the matrix elements of $\{A_i\}_{i=1}^n$  and $\{a_{i}\}_{i=1}^{mn}$ ($A_0$ is always invariant), and define birational transformations of the varieties of $\{A_i\}$ with given $\{a_{i}\}$ and $\{d_j\}$. 

There exist remarkable subgroups $\Z^n\subset \Z^{m(n+1)-1}$ which define birational transformations on the space of all $A(z)$ (with fixed $A_0$ and with no restrictions on the roots of $\det A(z)$), but to see that we need to parameterize $A(z)$ differently. 

To define the new coordinates, we split the eigenvalues of $A(z)$ into $n$ groups of $m$ numbers each:
$$
\{a_1,\dots,a_{mn}\}=\{a_1^{(1)},\dots,a_m^{(1)}\}\cup\dots\cup\{a_1^{(n)},\dots,a_m^{(n)}\}.
$$
The splitting may be arbitrary.
Then we define $B_i$ to be the uniquely determined (remember, everything is generic) element of $\Mat$ with eigenvalues $\left\{a_j^{(i)}\right\}_{j=1}^m$, such that $z-B_i$ is a right divisor of $A(z)$:
$$
A(z)=(A_0z^{n-1}+A_1'z^{n-1}+\dots+A_{n-1}')(z-B_i).
$$
The matrix elements of $\{B_i\}_{i=1}^n$ are the new coordinates on the space of $A(z)$. 

The action of the subgroup $\Z^n$ mentioned above consists of shifting the eigenvalues in any group by the same integer assigned to this group, and also shifting the exponents $\{d_i\}$ by the same integer (which is equal to minus the sum of the group shifts). If we denote by $\{B_i(k_1,\dots,k_n)\}$ the result of applying $k\in\Z^n$ to $\{B_i\}$, then the following equations are satisfied:
$$
\gather
B_i(\dots)-B_i(\dots ,k_j+1,\dots)=
B_j(\dots)-B_j(\dots ,k_i+1,\dots),\tag 3\\
{B_j(\dots ,k_i+1,\dots)B_i(\dots)=B_i(\dots ,k_j+1,\dots)B_j(\dots)},\tag 4\\
B_i(k_1+1,\dots,k_n+1)=A_0^{-1}B_i(k_1,\dots,k_n)A_0-I,\tag 5
\endgather 
$$ 
where $i,j=1,\dots,n$, and dots in the arguments mean that other $k_l$'s remain unchanged. We call them the {\it difference Schlesinger equation} for the reasons that will be clarified below. Note that \tht{3} and \tht{4} can be rewritten as
$$
\bigl(z-B_i(\dots,k_j+1,\dots)\bigr)\bigl(z-B_j(\dots)\bigr)=
\bigl(z-B_j(\dots,k_i+1,\dots)\bigr)\bigl(z-B_i(\dots)\bigr).
$$

Independently of Birkhoff's general theory, we prove that the difference Schlesinger equations have a unique solution satisfying 
$$
Sp(B_i(k_1,\dots,k_n))=Sp(B_i)-k_i,\qquad i=1,\dots,n,
\tag 6
$$
for an arbitrary nondegenerate $A_0$ and generic initial conditions $\{B_i=B_i(0)\}$. (The notation means that the eigenvalues of $B_i(k)$ are
equal to those of $B_i$ shifted by $-k_i$.)
Moreover, the matrix elements of this solution are rational functions in the matrix elements of the initial conditions. This is our second result.

In order to prove this claim, we introduce yet another set of coordinates on $A(z)$ with fixed $A_0$, which is related to $\{B_i\}$ by a birational transformation. It consists of matrices $C_i\in\Mat$ with
$Sp(C_i)=Sp(B_i)$ such that
$$
A(z)=A_0(z-C_1)\cdots(z-C_n).
$$
In these coordinates, the action of $\Z^n$ is described by the relations
$$
\gathered
\bigl(z+1-C_i\bigr)
\cdots
\bigl(z+1-C_n\bigr)
A_0\bigl(z-C_1\bigr)\cdots
\bigl(z-C_{i-1}\bigr)
\\=
\bigl(z+1-\wt C_{i+1}\bigr)
\cdots
\bigl(z+1-\wt C_n\bigr) 
A_0\bigl(z-\wt C_1\bigr)\cdots
\bigl(z-\wt C_{i}\bigr),
\\
C_j=C_j(k_1,\dots,k_n),\quad \wt C_j=C_j(k_1,\dots,k_{i-1},k_i+1,k_{i+1},\dots,k_n)\text{  for all  }j.
\endgathered
\tag 7
$$
Again, we prove that there exists a unique solution to these equations satisfying
$Sp(C_i(k))=Sp(C_i)-k_i$, for an arbitrary invertible $A_0$ and generic $\{C_i=C_i(0)\}$. The solution is rational in the matrix elements of the initial conditions.

The difference Schlesinger equations have an autonomous limit which consists of \tht{3}, \tht{4}, and
$$
\gather
B_i(k_1+1,\dots,k_n+1)=A_0^{-1}B_i(k_1,\dots,k_n)A_0,\tag 5-aut\\
Sp(B_i(k_1,\dots,k_n))=Sp(B_i),\qquad i=1,\dots,n.
\tag 6-aut
\endgather
$$
The equation \tht{7} then turns into
$$
\gathered
\bigl(z-C_i\bigr)
\cdots
\bigl(z-C_n\bigr)
A_0\bigl(z-C_1\bigr)\cdots
\bigl(z-C_{i-1}\bigr)
\\=
\bigl(z-\wt C_{i+1}\bigr)
\cdots
\bigl(z-\wt C_n\bigr) 
A_0\bigl(z-\wt C_1\bigr)\cdots
\bigl(z-\wt C_{i}\bigr).
\endgathered
\tag 7-aut
$$

The solutions of these equations were essentially obtained in \cite{V} via a general construction of commuting flows associated with set-theoretical solutions of the quantum Yang-Baxter equation, see \cite{V} for details and references.

The autonomous equations can also be explicitly solved in terms of abelian functions associated with the spectral curve
$\{(z,w): \det(A(z)-wI)=0\}$,\footnote{It is easy to see that the curve is invariant under the flows.} very much in the spirit of \cite{MV, \S1.5}. We hope to explain the details in a separate publication.

The whole subject bears a strong similarity (and not just by name!) to the theory of isomonodromy deformations of linear systems of differential equations with rational coefficients:
$$
\frac{d\Y(\ze)}{d\ze}=\left(\B_{\infty}+\sum_{k=1}^n\frac{\B_i}{\ze-x_i}\right)\Y(\ze),
\tag 8
$$
which was developed by Schlesinger around 1912 and generalized by Jimbo, Miwa, and Ueno in \cite{JMU}, \cite{JM} to the case of higher order singularities. If we analytically continue any fixed (say, normalized at a given point) solution $\Y(\ze)$ of \tht{8} along a closed path $\gamma$ in $\C$ avoiding the singular points $\{x_k\}$ then the columns of $\Y$ will change into their linear combinations: $\Y\mapsto \Y M_\gamma$. Here $M_\gamma$ is a constant invertible matrix which depends only on the homotopy class of $\gamma$. It is called the {\it monodromy matrix} corresponding to $\gamma$. The monodromy matrices define a linear representation of the fundamental group of $\C$ with $n$ punctures. The basic isomonodromy problem is to change the differential equation \tht{8} so that the monodromy representation remains invariant. 

There exist isomonodromy deformations of two types: continuous ones, when $x_i$ move in the complex plane and $\B_i=\B_i(x)$ form a solution of a system of partial differential equations called {\it Schlesinger equations}, and discrete ones (called {\it Schlesinger transformations}), which shift the eigenvalues of $\B_i$ and exponents of 
$\Y(\ze)$ at $\ze=\infty$ by integers with the total sum of shifts equal to 0.

We prove that in the limit when 
$$
B_i=x_i\epsilon^{-1}+\Cal\B_i,\qquad \epsilon\to 0,
$$
our action of $\Z^{m(n+1)-1}$ in the discrete case converges to the action of Schlesinger transformations on $\B_i$. This is our third result.

Furthermore, we argue that the ``long-time'' asymptotics of the $\Z^n$-action in the discrete case  (that is, the asymptotics of $B_i([x_1\epsilon^{-1}],\dots,[x_n\epsilon^{-1}])$), $\epsilon$ small, is described by the corresponding solution of the Schlesinger equations. More exactly, we conjecture that the following is true.

 Take $B_i=B_i(\epsilon)\in\Mat$, $i=1,\dots,n$, such that
$$
B_i(\epsilon)-y_i\epsilon^{-1}+\B_i\to 0,\qquad \epsilon\to 0.
$$
Let $B_i(k_1,\dots,k_n)$ be the solution of the difference Schlesinger equations \tht{3.1}-\tht{3.3} with the initial conditions $\{B_i(0)=B_i\}$, and let $\B_i(x_1,\dots,x_n)$ be the solution of the classical Schlesinger equations \tht{5.4} with the initial conditions
$
\{\B_i(y_1,\dots,y_n)=\B_i\}.
$
Then for any $x_1,\dots,x_n\in\Bbb R$ and $i=1,\dots,n$, we have
$$
B_i\bigl([x_1\epsilon^{-1}],\dots,[x_n\epsilon^{-1}]\bigr)+
(x_i-y_i)\epsilon^{-1}+\B_i(y_1-x_1,\dots,y_n-x_n)\to 0,\qquad \epsilon\to 0.
$$

In support of this conjecture, we explicitly show that the difference Schlesinger equations converge to the conventional Schlesinger equations in the limit $\epsilon\to 0$.

Note that the monodromy representation of $\pi_1(\C\setminus\{x_1,\dots,x_n\})$ which provides the integrals of motion for the Schlesinger flows, has no obvious analog in the discrete situation. On the other hand, the obvious differential analog of the periodic matrix $P$, which contains all integrals of motion in the case of difference equations, gives only the monodromy information at infinity and does not carry any information about local monodromies around the poles $x_1,\dots,x_n$.

Most of the results of the present paper can be carried over to the case of $q$-difference equations of the form $Y(qz)=A(z)Y(z)$. The $q$-difference Schlesinger equations are, cf. \tht{3-6},
$$
\gather
B_i(\dots)-B_i(\dots ,q^{k_j+1},\dots)=
B_j(\dots)-B_j(\dots ,q^{k_i+1},\dots),\tag 3q\\
{B_j(\dots ,q^{k_i+1},\dots)B_i(\dots)=B_i(\dots ,q^{k_j+1},\dots)B_j(\dots)},\tag 4q\\
B_i(q^{k_1+1},\dots,q^{k_n+1})=q^{-1}A_0^{-1}B_i(q^{k_1},\dots,q^{k_n})A_0,\tag 5q\\
Sp(B_i(q^{k_1},\dots,q^{k_n}))=q^{-k_i}Sp(B_i),\qquad i=1,\dots,n.
\tag 6q
\endgather 
$$ 
The $q$-analog of \tht{7} takes the form
$$
\gather
\bigl(z-q^{-1}C_i\bigr)
\cdots
\bigl(z-q^{-1}C_n\bigr)
A_0\bigl(z-C_1\bigr)\cdots
\bigl(z-C_{i-1}\bigr)
\\=
\bigl(z-q^{-1}\wt C_{i+1}\bigr)
\cdots
\bigl(z-q^{-1}\wt C_n\bigr) 
A_0\bigl(z-\wt C_1\bigr)\cdots
\bigl(z-\wt C_{i}\bigr),\tag 7q
\\
C_j=C_j(q^{k_1},\dots,q^{k_n}),\quad \wt C_j=C_j(q^{k_1},\dots,q^{k_{i-1}},q^{k_i+1},q^{k_{i+1}},\dots,q^{k_n})\text{  for all  }j.
\endgather
$$
A more detailed exposition of the $q$-difference case will appear elsewhere.

Similarly to the classical case, see \cite{JM}, discrete Painlev\'e equations of \cite{JS}, \cite{Sak} can be obtained as reductions of the difference and $q$-difference Schlesinger equations when both $m$ (the size of matrices) and $n$ (the degree of the polynomial $A(z)$) are equal to two. For examples of such reductions see
\cite{Bor, \S3} for difference Painlev\'e II equation (dPII), \cite{Bor, \S6} and \cite{BB, \S9} for dPIV and dPV, and \cite{BB, \S10} for $q$-PVI. This subject still remains to be thoroughly studied.

As was mentioned before, the difference and $q$-difference Schlesinger equations can be used to compute the gap probabilities for certain probabilistic models. We conclude this Introduction by giving an example of such a model. We define the {\it Hahn orthogonal polynomial ensemble} as a probability measure on all $l$-point subsets of $\{0,1,\dots,N\}$, $N>l>0$, such that 
$$
\operatorname{Prob}\{(x_1,\dots,x_l)\}=\operatorname{const}\cdot
\prod_{1\le i<j\le l}(x_i-x_j)^2\cdot \prod_{i=1}^l w(x_i),
$$ 
where $w(x)$ is the weight function for the classical Hahn orthogonal polynomials:
$$
w(x)={{\alpha+x}\choose {x}}{\beta+N-x\choose N-x},\qquad \alpha,\beta>-1\text{  or  } \alpha,\beta<-N.
$$
This ensemble came up recently in harmonic analysis on the infinite--dimensional unitary group \cite{BO, \S11} and in a statistical description of tilings of a hexagon by rhombi \cite{Joh, \S4}. 

The quantity of interest is the probability that the point configuration $(x_1,\dots,x_l)$ does not intersect a disjoint union of intervals $[k_1,k_2]\sqcup\dots\sqcup[k_{2s-1},k_{2s}]$. As a function in the endpoints $k_1,\dots,k_{2s}\in\{0,1,\dots,N\}$, this probability can be expressed through a solution of the difference Schlesinger equations \tht{3}-\tht{6} for $2\times 2$ matrices with $n=\deg A(z)=s+2$, $A_0=I$,
$$
\gathered
Sp(B_i)=\{-k_i,-k_i\}, \quad i=1,\dots,2s,\\ Sp(B_{2s+1})\sqcup Sp(B_{2s+2})=\{0,-\alpha,N+1,N+1+\beta\},
\endgathered
$$
and with certain explicit initial conditions. The equations are also suitable for numerical computations, and we refer to \cite{BB, \S12} for examples of those in the case of one interval gap. 

I am very grateful to P.~Deift, P.~Deligne, B.~Dubrovin, A.~Its, D.~Kazhdan, I.~Krichever, G.~Olshanski, V.~Retakh, and A.~Veselov for interesting and 
helpful discussions. 

This research was partially conducted during the period the author served as a Clay Mathematics Institute Long-Term Prize Fellow.

\head
1. Birkhoff's theory 
\endhead
Consider a matrix linear difference equation of the first order
$$
Y(z+1)=A(z)Y(z).
\tag 1.1
$$
Here $A:\C\to\Mat$ is a rational function (i.e., all matrix elements of $A(z)$ are rational functions of $z$) and $m\ge 1$. We are interested in matrix meromorphic solutions $Y:\C\to\Mat$ of this equation.

Let $n$ be the order of the pole of $A(z)$ at infinity, that is,
$$
A(z)=A_0z^n+A_1z^{n-1}+\text{  lower order terms }.
$$

We assume that \tht{1.1} has a formal solution of the form
$$
Y(z)=z^{nz}e^{-nz}\left(\hat Y_0+\frac{\hat Y_1}z+\frac{\hat Y_2}{z^2}+\dots\right)
\diag\left(\,\rho_1^z\,z^{d_1},\dots,\rho_m^z\,z^{d_m}\right)	 
\tag 1.2
$$
with $\rho_1,\dots,\rho_m\ne 0$ and $\det \hat Y_0\ne 0$.\footnote{Substituting \tht{1.2} in \tht{1.1} we use the expansion
$\left(\frac{z+1}z\right)^{nz}=e^{nz\ln(1+z^{-1})}=e^n-\frac{ne^n}{2z}+\dots$ to compare the two sides.}

It is easy to see that if such a formal solution exists then $\rho_1,\dots,\rho_m$ must be the eigenvalues of $A_0$, and the columns of $\hat Y_0$ must be the corresponding eigenvectors of $A_0$. 

Note that for any invertible $T\in \Mat$, $(TY)(z)$ solves the equation 
$$
(TY)(z+1)=(TA(z)T^{-1})\,(TY)(z).
$$
Thus, if $A_0$ is diagonalizable, we may assume that it is diagonal without loss of generality. Similarly, if $A_0=I$ and $A_1$ is diagonalizable, we may assume that $A_1$ is diagonal.

\proclaim{Proposition 1.1} If $A_0=\diag(\rho_1,\dots,\rho_m)$, where $\{\rho_i\}_{i=1}^m$ are nonzero and pairwise distinct, then there exists a unique formal solution of \tht{1.1} of the form \tht{1.2} with $\hat Y_0=I$. 
\endproclaim
\demo{Proof} It suffice to consider the case $n=0$; the general case is reduced to it by considering $(\Gamma(z))^nY(z)$ instead of $Y(z)$, because $$
\Gamma(z)=\sqrt{2\pi}\,z^{z-\frac 12}e^{-z}\left(1+\frac 1{12}\,z^{-1}+\dots\right).
$$
(More precisely, this expression formally solves $\Gamma(z+1)=z\Gamma(z)$.)

Thus, we assume $n=0$. Then we substitute \tht{1.2} into \tht{1.1} and compute $\hat Y_k$ one by one by equating the coefficients of $z^{-l}$, $l=0,1,\dots$\,. If $\hat Y_0=I$ then the constant coefficients of both sides are trivially equal. The coefficients of $z^{-1}$ give
$$
\hat Y_1 A_0+\diag(\rho_1d_1,\dots,\rho_md_m)=A_0\hat Y_1+A_1.
\tag 1.3
$$
This equality uniquely determines $\{d_i\}$ and the off-diagonal entries of $\hat Y_1$, because 
$$
{[\hat Y_1,A_0]}_{ij}=(\rho_j-\rho_i){(\hat Y_1)}_{ij}\,.
$$
Comparing the coefficients of $z^{-2}$ we obtain
$$
(\hat Y_2-\hat Y_1) A_0+\hat Y_1\diag(\rho_1d_1,\dots,\rho_md_m)+\ldots=
A_0\hat Y_2+A_1\hat Y_1+\dots\,,
$$
where the dots stand for the terms which we already know (that is, those which depend only on $\rho_i$'s, $d_i$'s, $A_i$'s, and $\hat Y_0=I$). Since the diagonal values of $A_1$ are exactly $\rho_1d_1,\dots\rho_nd_n$ by \tht{1.3}, we see that we can uniquely determine the diagonal elements of $\hat Y_1$ and the off-diagonal elements of $\hat Y_2$ from the last equality.

Now let us assume that we already determined $\hat Y_1,\dots,\hat Y_{l-2}$ and the off-diagonal entries of $\hat Y_{l-1}$ by satisfying \tht{1.1} up to order $l-1$. Then comparing the coefficients of $z^{-l}$ we obtain
$$
(\hat Y_{l}-(l-1)\hat Y_{l-1})A_0+
\hat Y_{l-1}\diag(\rho_1d_1,\dots,\rho_md_m)+\ldots=A_0\hat Y_l+A_1\hat Y_{l-1}+\dots\,,
$$
where the dots denote the terms depending only on $\rho_i$'s, $d_i$'s, $A_i$'s, and $\hat Y_0,\dots ,\hat Y_{l-2}$. This equality allows us to compute the diagonal entries of $Y_{l-1}$ and the off-diagonal entries of $Y_l$. Induction on $l$ completes the proof.\qed
\enddemo

The condition that the eigenvalues of $A_0$ are distinct is not necessary for the existence of the asymptotic solution, as our next proposition shows.

\proclaim{Proposition 1.2} Assume that $A_0=I$ and $A_1=\diag(r_1,\dots,r_n)$ where $r_i-r_j\notin \{\pm1,\pm2,\dots\}$ for all $i,j=1,\dots,n$. Then there exists a unique formal solution of \tht{1.1} of the form \tht{1.2} with $\hat Y_0=I$. 
\endproclaim
\demo{Proof} As in the proof of Proposition 1.1, we may assume that $n=0$. Comparing constant coefficients we see that $\rho_1=\dots=\rho_m=1$.
Then equating the coefficients of $z^{-1}$ we find that $d_i=r_i$, $i=1,\dots,m$. Furthermore, equating the coefficients of $z^{-l}$, $l\ge 2$ we find that
$$
[\hat Y_{l-1},A_1]-(l-1)\hat Y_{l-1}
$$
is expressible in terms of $A_i$'s and $\hat Y_1,\dots,\hat Y_{l-2}$.
This allows us to compute all $\hat Y_i$'s recursively.\qed
\enddemo

We call two complex numbers $z_1$ and $z_2$ {\it congruent} if $z_1-z_2\in\Bbb Z$.

\proclaim{Theorem 1.3 (G.~D.~Birkhoff \cite{Bi1, Theorem III})}
Assume that 
$$
\gather
A_0=\diag(\rho_1,\dots,\rho_m),\\
\rho_i\ne 0,\quad i=1,\dots,m,\qquad \rho_i/\rho_j\notin\Bbb R \ \text{  for all  }\  i\ne j.
\endgather
$$
Then there exist unique solutions $Y^l(z)$ 
{\rm(}$Y^r(z)${\rm )} of \tht{1.1} such that 

(a) the function $Y^l(z)$ 
{\rm(}$Y^r(z)${\rm )} is analytic throughout the complex plane except possibly for poles to the right {\rm (}left\,{\rm )} of and congruent to the poles of
$A(z)$ {\rm(}respectively, $A^{-1}(z-1)${\rm )};

(b) in any left {\rm (}right\,{\rm )} half plane $Y^l(z)$ {\rm(}$Y^r(z)${\rm )} is asymptotically represented by the right--hand side of \tht{1.2}.
\endproclaim

\example{Remark 1.4}
Part (b) of the theorem means that for any $k=0,1,\dots$,
$$
\left|Y^{l,r}(z)\,z^{-nz}e^{nz}
\diag(\rho_1^{-z}z^{-d_1},\dots,\rho_m^{-z}z^{-d_m})
-\hat Y_0-\frac{\hat Y_1}z-\dots-\frac{\hat Y_{k-1}}{z^{k-1}}\right|\le \frac {const}{z^{k}}
$$
for large $|z|$ in the corresponding domain.

Theorem holds for any (fixed) choices 
of branches of $\ln(z)$ in the left and right half planes for evaluating
$z^{-nz}=e^{-nz\ln(z)}$ and $z^{-d_k}=e^{-d_k\ln(z)}$,
and of a branch of $\ln(\rho)$ with a cut not passing through $\rho_1,\dots,\rho_m$ for evaluating $\rho_k^{-z}=e^{-z\ln\rho_k}$. Changing these branches yields the multiplication of $Y^{l,r}(z)$ by a diagonal periodic matrix on the right.
\endexample

\example{Remark 1.5} Birkhoff states Theorem 1.3 under a more general assumption: he only assumes that the equation \tht{1.1} has a formal solution of the form \tht{1.2}. However, as pointed out by P.~Deligne, Birkhoff's proof has a flaw in case one of the ratios $\rho_i/\rho_j$ is real. The following counterexample was
kindly communicated to me by Prof. Deligne.

Consider the equation \tht{1.1} with $m=2$ and
$$
A(z)=\bmatrix 1&1/z\\ 0& 1/e\endbmatrix.
$$
The formal solution \tht{1.2} has the form
$$
Y(z)=\left(I+\bmatrix 0&a\\0&0\endbmatrix z^{-1}+\dots\right)\bmatrix 1&0\\0&e^{-z}\endbmatrix
$$
with $a=e/(1-e)$.

Actual solutions that we care about have the form
$$
Y(z)=\bmatrix 1&u(z)\\ 0&e^{-z}\endbmatrix
$$
where $u(z)$ is a solution of 
$u(z+1)=u(z)+e^{-z}/z$. In a right half-plane we can take
$$
u^r(z)=-\sum_{n=0}^\infty\frac{e^{-(z+n)}}{z+n}\,.
$$
The first order approximation of $u^r(z)$ anywhere except near nonpositive integers is 
$$
u^r(z)\sim -\sum_{n=0}^\infty \frac{e^{-(z+n)}}z=\frac {ae^{-z}}z\,.
$$
Next terms can be obtained by expanding $1/(z+n)$.

In order to obtain a solution which behaves well on the left, it suffices to cancel the poles:
$$
u^l(z)=u^r(z)+\frac {2\pi i}{e^{2\pi i z}-1}\,.
$$
The corresponding solution $Y^l(z)$ has needed asymptotics in sectors of the form $\pi/2+\epsilon<\arg z<3\pi/2+\epsilon$, but it has wrong asymptotic behavior as
$z\to +i\infty$. Indeed, $\lim_{z\to +i\infty} u^l(z)=-2\pi i$.

On the other hand, we can take 
$$
\wt u^l(z)=u^l(z)+2\pi i=u^r(z)+\frac {2\pi i\,e^{2\pi iz}}{e^{2\pi i z}-1}\,,
$$
which has the correct asymptotic behavior in $\pi/2-\epsilon<\arg z<3\pi/2-\epsilon$, but fails to have the needed asymptotics at $-i\infty$.
\endexample

\example{Remark 1.6} In the case when $|\rho_1|>|\rho_2|>\dots>|\rho_m|>0$, a result similar to Theorem 1.3 was independently proved by R.~D.~Carmichael \cite{C}. He considered the asymptotics of solutions along lines parallel to the real axis only. Birkhoff also referred to \cite{N} and \cite{G} where similar results had been proved somewhat earlier.
\endexample

Now let us restrict ourselves to the case when $A(z)$ is a polynomial in $z$. The general case of rational $A(z)$ is reduced to the polynomial case by the following transformation. If $(z-x_1)\cdots(z-x_s)$ is the
common denominator of $\{A_{kl}(z)\}$ (the matrix elements of $A(z)$), then 
$$
\bar Y(z)=\Gamma(z-x_1)\cdots\Gamma(z-x_s)\cdot Y(z)
$$
solves $\bar Y(z+1)=\bar A(z)\bar Y(z)$ with polynomial
$$
\bar A(z)=(z-x_1)\cdots (z-x_s) A(z).
$$

Note that the ratio $P(z)=(Y^r(z))^{-1}\,Y^l(z)$ is a periodic function. (The relation $P(z+1)=P(z)$ immediately follows from the fact that $Y^{l,r}$ solve \tht{1.1}.) From now on let us fix the branches of $\ln(z)$ in the left and right half planes mentioned in Remark 1.4 so that they coincide in the upper half plane. Then the structure of $P(z)$ can be described more precisely.

\proclaim{Theorem 1.7 (G.~D.~Birkhoff \cite{Bi1, Theorem IV})}
In the assumptions of Theorem 1.3, the matrix elements $p_{kl}(z)$ of the periodic matrix $P(z)=(Y^r(z))^{-1}Y^l(z)$ have the form
$$
\gathered
p_{kk}(z)=1+c_{kk}^{(1)}e^{2\pi iz}+\dots+c_{kk}^{(n-1)}e^{2\pi(n-1)iz}+
e^{2\pi i d_k}e^{2\pi n iz},\\
p_{kl}(z)=e^{2\pi \lambda_{kl}z}
\left(c_{kl}^{(0)}+c_{kl}^{(1)}e^{2\pi iz}+\dots+c_{kl}^{(n-1)}e^{2\pi(n-1)iz}\right) \quad (k\ne l),
\endgathered
$$
where $c_{kl}^{(s)}$ are some constants, and $\lambda_{kl}$ denotes the least integer as great as the real part of $(\ln(\rho_l)-\ln(\rho_k))/2\pi i$.
\endproclaim

Thus, starting with a matrix polynomial $A(z)=A_0z^n+A_1z^{n-1}+\dots+A_n$ with nondegenerate $A_0=\diag(\rho_1,\dots,\rho_m)$, $\rho_k\ne\rho_l$ for $k\ne l$, we construct the {\it characteristic constants} $\{d_k\}$, $\{c_{kl}^{(s)}\}$ using Proposition 1.1 and Theorems 1.3, 1.7.

Note that the total number of characteristic constants is exactly the same as the number of matrix elements in matrices $A_1,\dots,A_n$. Thus, it is natural to ask whether the map 
$$
(A_1,\dots,A_n)\mapsto \left(\{d_k\},\, \{c_{kl}^{(s)}\}\right)
$$
is injective or surjective (the constants $\rho_1,\dots,\rho_n$ are being fixed). The following partial results are available.

\proclaim{Theorem 1.8 (G.~D.~Birkhoff \cite{Bi2, \S17})}  For any nonzero $\rho_1,\dots,\rho_m$, $\rho_i/\rho_j\notin\Bbb R$ for $i\ne j$, there exist matrices $A_1,\dots,A_n$ such that the equation \tht{1.1} with $A_0=\diag(\rho_1,\dots\rho_m)$ either possesses prescribed characteristic  
constants $\{d_k\}$, $\{c_{kl}^{(s)}\}$, or else constants
$\{d_k+l_k\},\, \{c_{kl}^{(s)}\}$, where $l_1,\dots,l_m$ are integers.
\endproclaim

\proclaim{Theorem 1.9 (G.~D.~Birkhoff \cite{Bi1, Theorem VII})}
Assume that we have two matrix polynomials $A'(z)=A_0'z^n+\dots+A_n'$ and $A''(z)=A_0''z^n+\dots+A_n''$ with
$$
A_0'=A_0''=\diag(\rho_1,\dots\rho_m), \qquad \rho_k\ne0,\quad \rho_k/\rho_l\notin \Bbb R\text{   for   }k\ne l,
$$
such that the sets of the characteristic constants for the equations
$Y'(z+1)=A'(z)Y'(z)$ and $Y''(z+1)=A''(z)Y''(z)$ are equal. Then there exists a rational matrix $R(z)$ such that 
$$
A''(z)=R(z+1)A'(z)R^{-1}(z),
$$
and the left and right canonical solutions $Y^{l,r}$ of the second equation can be obtained from those of the first equation by multiplication by $R$ on the left:
$$
(Y'')^{l,r}=R\,(Y')^{l,r}.
$$
\endproclaim

\head 2. Isomonodromy transformations
\endhead

The goal of this section is to construct explicitly, for given $A(z)$, rational matrices $R(z)$ such that the transformation $A(z)\mapsto R(z+1)A(z)R^{-1}(z)$, cf. Theorem 1.9 above, preserves the characteristic constants (more generally, preserves $\{c_{kl}^{(s)}\}$ and shifts $d_k$'s by integers).

Let $A(z)$ be a matrix polynomial of degree $n\ge 1$, $A_0=\diag(\rho_1,\dots,\rho_m)$, and $\rho_i$'s are nonzero and their ratios are not real. Fix $mn$ complex numbers $a_1,\dots,a_{mn}$ such that $a_i-a_j\notin \Z$ for any $i\ne j$. 
Denote by $\M(a_1,\dots,a_{mn}; d_1,\dots,d_m)$ the algebraic variety of all $n$-tuples of $m$ by $m$ matrices $A_1,\dots,A_n$ such that the scalar polynomial
$$
\det A(z)=\det(A_0z^n+A_1z^{n-1}+\dots+A_n)
$$
of degree $mn$ has roots $a_1,\dots,a_{mn}$, and $\rho_i\left(d_i-\frac n2\right)={(A_1)}_{ii}$ (this comes from the analog of \tht{1.3} for arbitrary $n$). 

\proclaim{Theorem 2.1} For any $\kappa_1,\dots,\kappa_{mn}\in \Z$,
$\delta_1,\dots,\delta_m\in \Z$,
$$
\sum_{i=1}^{mn} \kappa_i+\sum_{j=1}^{m}\delta_j=0,
$$ 
there exists a nonempty Zariski open subset $\A$ of $\M(a_1,\dots,a_{mn};d_1,\dots,d_m)$ such that for any $(A_1,\dots,A_n)\in\A$ there exists a unique rational matrix $R(z)$ with the following properties:
$$
\gather
\wt A(z)=R(z+1)A(z)R^{-1}(z)=\wt A_0 z^n+\wt A_1z^{n-1}+\dots+\wt A_n,\quad \wt A_0=A_0, 
\\
(\wt A_1,\dots,\wt A_n)\in \M(a_1+\kappa_1,\dots,a_{mn}+\kappa_{mn};d_1+\delta_1,\dots, d_m+\delta_m),
\endgather
$$
and the left and right canonical solutions of $\wt Y(z+1)=\wt A(z)\wt Y(z)$ have the form
$$
\wt Y^{l,r}=R\, Y^{l,r},
$$
where $Y^{l,r}$ are left and right canonical solutions of $Y(z+1)=A(z)Y(z)$.

The map $(A_1,\dots,A_n)\mapsto (\wt A_1,\dots,\wt A_n)$ is a birational map of algebraic varieties.
\endproclaim

\example{Remark 2.2}  The theorem implies that the characteristic constants $\{c_{kl}^{(s)}\}$ for the difference equations with coefficients $A$ and $\wt A$ are
the same, while the constants $d_k$ are being shifted by $\delta_k\in\Z$.

Note also that if we require that all $d_k$'s do not change then, by virtue of Theorem 1.9, Theorem 2.1 provides all possible transformations which preserve  the characteristic constants. Indeed, if $A''(z)=R(z+1)A'(z)R^{-1}(z)$ then zeros of $\det A''(z)$ must be equal to those of $\det A'(z)$ shifted by integers.
\endexample
\demo{Proof} Let us proof the uniqueness of $R$ first. Assume that there exist two rational matrices $R_1$ and $R_2$ with needed properties. 
This means, in particular, that the determinants of the matrices
$$
\wt A^{(1)}=R_1(z+1)A(z)R^{-1}_1(z)\quad\text{  and  }\quad
\wt A^{(2)}=R_2(z+1)A(z)R^{-1}_2(z)
$$
vanish at the same set of $mn$ points $\wt a_{i}=a_i+\kappa_i$, none of which are different by an integer. Denote by $\wt Y^r_1=R_1Y^r$ and $\wt Y_2^r=R_2Y^r$ the right canonical solutions of the corresponding equations. Then
$\wt Y_1^r(\wt Y_2^r)^{-1}=R_1R_2^{-1}$ is a rational matrix which tends to $I$ at infinity. Moreover,
$$
\left(\wt Y_1^r(\wt Y_2^r)^{-1}\right)
(z+1)=\wt A^{(1)}(z)\, \left(\wt Y_1^r(\wt Y_2^r)^{-1}\right)(z)
\,\left(\wt A^{(2)}(z)\right)^{-1}.
$$
Since $\wt Y_1^r(\wt Y_2^r)^{-1}$ is holomorphic for $\Re z\ll 0$, the equation above implies that this function may only have poles at the points which are congruent to $\wt a_i$ (zeros of $\det \wt A^{(2)}(z)$) and to the right of them. (Recall that two complex numbers are congruent if their difference is an integer.) But since $\wt Y_1^r(\wt Y_2^r)^{-1}$ is also holomorphic for $\Re z\gg 0$, the same equation rewritten as
$$
\left(\wt Y_1^r(\wt Y_2^r)^{-1}\right)
(z)=\left(\wt A^{(1)}(z)\right)^{-1}\, \left(\wt Y_1^r(\wt Y_2^r)^{-1}\right)(z+1)
\,\wt A^{(2)}(z)
$$
implies that this function may only have poles at the points $\wt a_i$ (zeros of $\det\wt A^{(1)}(z)$) or at the points congruent to them and to the left of them. Thus, $\wt Y_1^r(\wt Y_2^r)^{-1}=R_1R_2^{-1}$ is entire, and by Liouville's theorem it is identically equal to $I$. The proof of uniqueness is complete.

To prove the existence we note, first of all, that it suffices to provide
a proof if one of $\kappa_i$'s is equal to $\pm 1$ and one of $\delta_j$'s is equal to $\mp 1$ with all other $\kappa$'s and $\delta$'s being equal to zero. The proof will consist of several steps.

\proclaim{Lemma 2.3} Let $A(z)$ be an $m$ by $m$ matrix--valued function holomorphic near $z=a$, and $\det A(z)=c(z-a)+O\left((z-a)^2\right)$ as $z\to a$, where $c\ne 0$. Then there exists a unique (up to a constant) nonzero vector $v\in \C^m$ such that $A(a)v=0$. 
Furthermore, if $B(z)$ is another matrix--valued function which is holomorphic near $z=a$, then $(BA^{-1})(z)$ is holomorphic near $z=a$ if and only if
$B(a)v=0$.
\endproclaim
\demo{Proof} Let us denote by $E_1$ the matrix unit which has 1 as its $(1,1)$-entry and 0 as all other entries. Since $\det A(a)=0$, there exists a nondegenerate constant matrix $C$ such that $A(a)CE_1=0$ (the first column of $C$ must be a 0-eigenvector of $A(a)$). This implies that
$$
H(z)=A(z)C(E_1(z-a)^{-1}+I-E_1)
$$
is holomorphic near $z=a$. On the other hand, $\det H(a)=c\det C\ne 0$.
Thus, $A(a)=H(a)(I-E_1)C^{-1}$ annihilate a vector $v$ if and only if 
$C^{-1}v$ is proportional to $(1,0,\dots,0)^t$. Hence, $v$ must be proportional to the first column of $C$. The proof of the first part of the lemma is complete.

To prove the second part, we notice that
$$
(BA^{-1})(z)=B(z)C(E_1(z-a)^{-1}+I-E_1)H^{-1}(z)
$$
which is bounded at $z=a$ if and only if $B(a)CE_1=0$. \qed
\enddemo

More generally, we will denote by $E_i$ the matrix unit defined by
$$
{(E_i)}_{kl}=\cases 1, &k=l=i,\\
                    0, &\text{otherwise}.
             \endcases
$$

\proclaim{Lemma 2.4 (\cite{JM, \S2 and Appendix A})} For any nonzero vector $v=(v_1,\dots,v_m)^t$, $Q\in\Mat$, $a\in \C$, and $i\in\{1,\dots,m\}$, there exists a linear matrix--valued function $R(z)=R_{-1}(z-a)+R_0$ with the properties
$$
\gather
R(z)\left(I+Qz^{-1}+O\left(z^{-2}\right)\right)z^{-E_i}
=I+O(z^{-1}),\quad z\to\infty,\\
R_0\,v=0,
\endgather
$$
if and only if $v_i\ne 0$. In this case, $R^{\pm 1}(z)$ is given by
$$
\gather
R(z)=E_i(z-a)+R_0,\quad R^{-1}(z)=I-E_i+R_1(z-a)^{-1},\quad \det R(z)=z-a,\\
{(R_0)}_{kl}=\cases v_i^{-1}\sum_{s\ne i}Q_{is}v_s,&k=l=i,\\
                    -Q_{il},& k=i,\,l\ne i,\\
		    -v_i^{-1}v_k,&k\ne i,\,l=i,\\
		    \delta_{kl},&k\ne i,\,l\ne i,
	     \endcases
\\
R_1=v_i^{-1}v\bmatrix Q_{i1},\dots,Q_{i,i-1},
 1 ,Q_{i,i+1},\dots,Q_{im}\endbmatrix.	     
\endgather
$$	         
\endproclaim	         
The proof is straightforward.

Now we return to the proof of Theorem 2.1. 

Assume that $\kappa_1=-1$, $\delta_i=1$ for some $i=1,\dots,m$, and all other $\kappa$'s and $\delta$'s are zero. Since $a_1$ is a simple root of $\det A(z)$, by Lemma 2.3 there exists a unique (up to a constant) vector $v$ such that
$A(a)v=0$. Clearly, the condition $v_i\ne 0$ defines a nonempty Zariski open subset of $\M(a_1,\dots,a_{mn};\delta_1,\dots;\delta_m)$. On this subset, let us take $R(z)$ to be the matrix afforded by Lemma 2.4 with $a=a_1$ and
$Q=\hat Y_1$ (we assume that $\hat Y_0=I$, see Proposition 1.1). Then
by the second part of Lemma 2.3, $(A(z)R^{-1}(z))^{-1}=R(z)A^{-1}(z)$ is holomorphic and invertible near $z=a_1$ (the invertibility follows from the fact that $\det R(z)A^{-1}(z)$ tends to a nonzero value as $z\to a_1$). Thus, 
$\wt A(z)=R(z+1)A(z)R^{-1}(z)$ is entire, hence, it is a polynomial. Since 
$$
\det\wt A(z)=\frac{z+1-a_1}{z-a_1}\,\det A(z)=c\,(z+1-a_1)(z-a_2)\cdots(z-a_{mn}), \quad c\ne 0,
$$
the degree of $A(z)$ is $\ge n$. Looking at the asymptotics at infinity, we see that $\deg A(z)\le n$, which means that $\wt A$ is a polynomial of degree $n$:
$$
\wt A(z) =\wt A_0z^n+\dots+\wt A_n, \quad \wt A_0\ne 0.
$$

Denote by $Y^{l,r}$ the left and right canonical solutions of $Y(z+1)=A(z)Y(z)$ (see Theorem 1.3 above). Then $\wt Y^{l,r}:=R\,Y^{l,r}$
are solutions of $\wt Y(z+1)=\wt A(z)\wt Y(z)$. Moreover, their asymptotics at infinity at any left (right) half plane, by Lemma 2.4, is given by an expansion of the form \tht{1.2} with $\hat {\wt Y_0}=I$,
$\wt\rho_k=\rho_k$ for all $k=1,\dots,m$, and 
$$
\wt d_k=\cases d_k+1,& k=i,\\
               d_k,&k\ne i.
	\endcases              
$$
This implies that $\wt A_0=\diag(\rho_1,\dots,\rho_m)$, and that $\wt Y^{l,r}$ are the left and right canonical solutions of the equation
$\wt Y(z+1)=\wt A(z)\wt Y(z)$. Indeed, their asymptotic expansion at infinity must also be a formal solution of the equation,
the fact that $\wt Y^{l,r}$ are holomorphic for $\Re z\ll 0$ ($\gg 0$) follows from the analogous property for $Y^{l,r}$, and the location of possible poles of $\wt Y^r$ is easily determined from the equation.

For a future reference let us also find a (unique up to a constant) vector $\wt v$ such that $\wt A^t(a_1-1)\,\wt v=0$. This means that $R^{-t}(a_1-1)A^t(a_1-1)R^t(a_1)\,\wt v=0$. Lemma 2.4 then implies that 
$$
\wt v=\bmatrix  {(\hat Y_1)}_{i1},\dots, {(\hat Y_1)}_{i,i-1},1,{(\hat Y_1)}_{i,i+1},\dots,{(\hat Y_1)}_{im}\endbmatrix^t
$$
is a solution. Note that $\wt v_i\ne 0$.

Now let us assume that $\kappa_1=1$ and $\delta_i=-1$ for some $i=1,\dots,m$. By Lemma 2.3, there exists a unique (up a to a constant) vector $w$ such that $A^t(a_1)w=0$. The condition $w_i\ne 0$ defines a nonempty Zariski open subset of $\M(a_1,\dots,a_{mn};\delta_1,\dots\delta_m)$.
On this subset, denote by $R'(z)$ the rational matrix-valued function afforded by Lemma 2.4 with $a=a_1$, $v=w$, and $Q=-\hat Y_1^t$ (again, we assume that $\hat Y_0=I$). Set 
$$
R(z):={(R')}^{-t}(z-1).
$$
Then by Lemma 2.4
$$
R(z)\left(I+\hat Y_1z^{-1}+O\left(z^{-2}\right)\right)z^{E_i}
=I+O(z^{-1}),\quad z\to\infty.
$$
Furthermore, by Lemma 2.3, $R^{-t}(z+1)A^{-t}(z)$ is holomorphic and invertible near $z=a_1$. Hence, $\wt A(z)=R(z+1)A(z)R^{-1}(z)$ is entire (note that $R^{-1}(z)=(R')^t(z-1)$ is linear in $z$). The rest of the argument is similar to the case $\kappa_1=-1$, $\delta_i=1$ considered above. 

Finding a solution $\wt w$ to $\wt A(a_1+1) \wt w=0$ is equivalent to finding a solution to $R'(a_1)\wt w=0$. One such solution has the form
$$
\wt w=\bmatrix  {-(\hat Y_1)}_{1i},\dots, -{(\hat Y_1)}_{i-1,i},1,-{(\hat Y_1)}_{i+1,i},\dots,-{(\hat Y_1)}_{mi}\endbmatrix^t
$$
and all others are proportional to it. Note that its $i$th coordinate is nonzero.

From what was said above, it is obvious that the image of the map 
$$
\M(a_1,\dots,a_{mn};\delta_1,\dots,\delta_m)\to
\M(a_1-1,\dots,a_{mn};\delta_1,\dots,\delta_i+1,\dots,\delta_m)
$$
is in the domain of definition of the map
$$
\M(a_1-1,\dots,a_{mn};\delta_1,\dots,\delta_i+1,\dots,\delta_m)\to
\M(a_1,\dots,a_{mn};\delta_1,\dots,\delta_m)
$$
and the other way around. On the other hand, the composition of these maps in either order must be equal to the identity map due to the uniqueness argument in the beginning of the proof. Hence, these maps are inverse to each other, and they establish a bijection between their domains of definition. The rationality of the maps follows from the explicit formula for $R(z)$ in Lemma 2.4. The proof of Theorem 2.1 is complete. \qed 
\enddemo

\example{Remark 2.5} Quite similarly to Lemma 2.4, the multiplier $R(z)$ can be computed in the cases when two $\kappa$'s are equal to $\pm 1$ or two $\delta$'s are equal to $\pm 1$ with all other $\kappa$'s and $\delta$'s being zero, cf. \cite{JM}.

Assume $\kappa_i=-1$ and $\kappa_j=1$. Denote by $v$ and $w$ the solutions of $A(a_i)\,v=0$ and $A^t(a_j)\, w=0$. Then $R$ exists if and only if $(v,w):=v^tw=w^tv\ne 0$, in which case 
$$
\gather
R(z)=I+\frac{R_0}{z-a_j-1}\,,\quad R^{-1}(z)=I-\frac{R_0}{z-a_i}\,,
\quad \det R(z)=\frac{z-a_i}{z-a_j-1}\,,\\
R_0=\frac{a_j-a_i+1}{(v,w)}\,vw^t.
\endgather
$$

Now assume $\delta_i=1$, $\delta_j=-1$. Then we must have
$\det R(z)=1$ and
$$
R(z)\left(I+\hat Y_1z^{-1}+\hat Y_2z^{-2}+O(z^{-3})\right)z^{E_j-E_i}=
I+O(z^{-1}),\quad z\to\infty.
$$
The solution exists if and only if ${(\hat Y_1)}_{ij}\ne 0$, in which case it has the form
$$
R(z)=E_iz+R_0,\quad R^{-1}(z)=E_jz+R_0^{-1},
$$
with ${(R_0)}_{kl}$ given by
$$
\matrix &l=i&l=j&l\ne i,j\\
                     k=i&\frac{-{(\hat Y_2)}_{ij}+\sum_{s\ne i}{(\hat Y_1)}_{is}{(\hat Y_1)}_{sj}}{{(\hat Y_1)}_{ij}}&-{(\hat Y_1)}_{ij}&
-{(\hat Y_1)}_{il}\\
k=j&\frac 1{{(\hat Y_1)}_{ij}}&0&0\\
k\ne i,j& -\frac{{(\hat Y_1)}_{kj}}{{(\hat Y_1)}_{ij}}&0&\delta_{kl}
\endmatrix		     
$$
and ${(R_0^{-1})}_{kl}$ given by
$$
\matrix &l=i&l=j&l\ne i,j\\
                     k=i&0&{(\hat Y_1)}_{ij}&0\\
k=j&-\frac 1{{(\hat Y_1)}_{ij}}&-\frac{{(\hat Y_2)}_{ij}}{{(\hat Y_1)}_{ij}}+{(\hat Y_1)}_{jj}& -\frac{{(\hat Y_1)}_{il}}{{(\hat Y_1)}_{ij}}\\
k\ne i,j& 0&{(\hat Y_1)}_{kj}&\delta_{kl}
\endmatrix
$$
\endexample

\head 3. Difference Schlesinger equations
\endhead

In this section we give a different description for the transformations $A\mapsto \wt A$ of Theorem 2.1 with 
$$
\kappa_{i_1}=\dots=\kappa_{i_m}=\pm 1,\quad
\delta_1=\dots=\delta_m=\mp 1,
$$
and all other $\kappa_i$'s are equal to zero, and for compositions of such transformations.

In what follows we always assume that 
our matrix polynomials $A(z)=A_0z^n+\dots$ have nondegenerate highest coefficients: $\det A_0\ne 0$. We also assume that $mn$ roots of the equation $\det A(z)=0$ are pairwise distinct; we will call them the {\it eigenvalues} of $A(z)$. For an eigenvalue $a$, there exists a (unique) nonzero vector $v$ such that $A(a)\,v=0$, see Lemma 2.3. We will call $v$ the {\it eigenvector} of $A(z)$ corresponding to the eigenvalue $a$.
The word {\it generic} everywhere below stands for ``belonging to a Zariski open subset'' of the corresponding algebraic variety.

We start with few simple preliminary lemmas. 

\proclaim{Lemma 3.1} The sets of eigenvalues and corresponding eigenvectors define $A(z)$ up to a multiplication by a constant nondegenerate matrix on the left.
\endproclaim
\demo{Proof} If there are two matrix polynomials $A'$ and $A''$ with same eigenvalues and eigenvectors, then $(A'(z))^{-1}A''(z)$ has no singularities in the finite plane. Moreover, since the degrees of $A'(z)$ and $A''(z)$ are equal, $(A'(z))^{-1}A''(z)\sim (A'_0)^{-1}A_0''$ as $z\to\infty$. Liouville's theorem concludes the proof.\qed
\enddemo

We will say that $z-B$, $B\in \Mat$, is a {\it right divisor} of $A(z)$ if $A(z)=\hat A(z)(z-B)$, where $\hat A(z)$ is a polynomial of degree $n-1$.

\proclaim{Lemma 3.2} A linear function $z-B$ is a right divisor of $A(z)$ if and only if 
$$
A_0B^n+A_1B^{n-1}+\dots+A_n=0.
$$
\endproclaim
\demo{Proof} See, e.g., \cite{GLR}.\qed
\enddemo

\proclaim{Lemma 3.3} Let $\alpha_1,\dots,\alpha_m$ be eigenvalues of $A(z)$ and $v_1,\dots,v_m$ be the corresponding eigenvectors. Assume that
$v_1,\dots,v_m$ are linearly independent. Take $B\in \Mat$ such that
$Bv_i=\alpha_iv_i$, $i=1,\dots,m$. Then $z-B$ is a right divisor of $A(z)$. Moreover, $B$ is uniquely defined by the conditions that $z-B$ is a right divisor of $A(z)$ and $Sp(B)=\{\alpha_1,\dots,\alpha_m\}$.
\endproclaim
\demo{Proof} For all $i=1,\dots,m$, we have
$$
(A_0B^n+A_1B^{n-1}+\dots+A_n)v_i=(A_0\alpha_i^n+A_1\alpha_i^{n-1}+\dots+A_n)v_i=A(\alpha_i)v_i=0.
$$
Lemma 3.2 shows that $z-B$ is a right divisor of $A(z)$.

To show uniqueness, assume that 
$$
A(z)=\hat A'(z)(z-B')=\hat A''(z)(z-B'').
$$
This implies $(A''(z))^{-1}A'(z)=(z-B'')(z-B')^{-1}$. Possible singularities of the right--hand side of this equality are $z=\alpha_i$, $i=1,\dots,m$, while possible singularities of the left--hand side are all other eigenvalues of $A(z)$. Since the eigenvalues of $A(z)$ are pairwise distinct, both sides are entire. But $(z-B'')(z-B')^{-1}$ tends to $I$ as $z\to\infty$. Hence, by Liouville's theorem, $B'=B''$. 
\qed
\enddemo

Now let us assume that the eigenvalues $a_1,\dots,a_{mn}$ of $A(z)$ are divided into $n$ groups of $m$ numbers:
$$
\{a_1,\dots,a_{mn}\}=\{a_1^{(1)},\dots,a_m^{(1)}\}\cup\dots\cup\{a_1^{(n)},\dots,a_m^{(n)}\}.
$$

Lemma 3.3 shows that for a generic $A(z)$ we can construct uniquely defined $B_1,\dots,B_n\in\Mat$ such that for any $i=1,\dots,n$, $Sp(B_i)=\{a_1^{(i)},\dots,a_m^{(i)}\}$ and $z-B_i$ is a right divisor of $A(z)$.\footnote{It is obvious that the condition on $A(z)$ used in Lemma 3.3, is an open condition. The corresponding set is nonempty because it contains diagonal $A(z)$ with $\bigl\{a_i^{(k)}\bigr\}$ being the roots of $A_{kk}(z)$. Similar remarks apply to all appearances of the word ``generic'' below.}
 By Lemma 3.1, $B_1,\dots,B_n$ define $A(z)$ uniquely up to a left constant factor, because the eigenvectors of $B_i$ must be eigenvectors of $A(z)$.

\proclaim{Lemma 3.4} For generic $B_1,\dots,B_n\in\Mat$ with
$Sp(B_i)=\{a_j^{(i)}\}$, there exists a unique monic degree $n$ polynomial $A(z)=z^n+A_1z^{n-1}+\dots$ such that $z-B_i$ are its right divisors.
The matrix elements of $A_1,\dots,A_n$ are rational functions of the matrix elements of $B_1,\dots,B_n$ and eigenvalues $\bigl\{a_j^{(i)}\bigr\}$.
\endproclaim
\example{Remark 3.5} 1. Later on we will show that, in fact, these rational functions do not depend on $\bigl\{a_j^{(i)}\bigr\}$.
\smallskip
\noindent 2. Clearly, the condition of $A(z)$ being monic can be replaced by the condition of $A(z)$ having a prescribed nondegenerate highest coefficient
$A_0$.
\endexample

\demo{Proof} The uniqueness follows from Lemma 3.1. To prove the existence part, we use induction on $n$. For $n=1$ the claim is obvious. Assume that we have already constructed $\hat A(z)=z^{n-1}+\hat A_1\dots$ such that
$B_1,\dots,B_{n-1}$ are its right divisors. Let $\{v_i\}$ be the eigenvectors of $B_n$ with eigenvalues $\{a^{(n)}_i\}$.
Set $w_i=\hat A(a^{(n)}_i)v_i$ and
take $X\in \Mat$ such that $X w_i=a^{(n)}_iw_i$ for all $i=1,\dots,m$.
(The vectors $\{w_i\}$ are linearly independent generically.) Then
$A(z)=(z-X)\hat A(z)$ has all needed properties. Indeed, we just need to check that $z-B_n$ is its right divisor (the rationality follows from the fact that computing the eigenvectors with known eigenvalues is a rational operation). For any $i=1,\dots,m$, we have
$$
\gathered
(B_n^n+A_1B_n^{n-1}+\dots+A_n)v_i=\left({(a_i^{(n)})}^n+A_1{(a_i^{(n)})}^{n-1}+\dots+A_n\right)v_i\\
=(a_i^{(n)}-X)\hat A(a_i^{(n)})v_i=
(a_i^{(n)}-X)w_i=0.
\endgathered
$$
Lemma 3.2 concludes the proof.\qed
\enddemo

Thus, we have a birational map between matrix polynomials $A(z)=A_0z^n+\dots$ with a fixed nondegenerate highest coefficient and
fixed mutually distinct eigenvalues divided into $n$ groups of $m$ numbers each, and sets of right divisors $\{z-B_1,\dots,z-B_n\}$ with $B_i$ having
the eigenvalues from the $i$th group. We will treat $\{B_i\}$ as a different set of coordinates for $A(z)$.

It turns out that in these coordinates some multipliers $R(z)$ of Theorem 2.1 take a very simple form. We will redenote by $\kappa^{(i)}_j$ the 
numbers $\kappa_1,\dots,\kappa_{mn}$ used in Theorem 2.1 in accordance with our new notation for the eigenvalues of $A(z)$. 
Denote the transformation of Theorem 2.1 with
$$
\kappa^{(i)}_j=-k_i\in \Z,\quad i=1,\dots,n,\,j=1,\dots,m;
\qquad \delta_1=\dots=\delta_m=\sum_{i=1}^n k_i
$$
by $\s(k_1,\dots,k_n)$.

\proclaim{Proposition 3.6} The multiplier $R(z)$ for $\s(0,\dots,0,\overset{(i)}\to 1,0,\dots,0)$ is equal to the right divisor $z-B_i$ of $A(z)$ corresponding to the eigenvalues $a^{(i)}_1,\dots,a^{(i)}_m$.
\endproclaim
\demo{Proof} It is easy to see that if $B_i$ has eigenvalues $a^{(i)}_1,\dots,a^{(i)}_m$, and $z-B_i$ is a right divisor of $A(z)$ then $R(z)=z-B_i$ satisfies all the conditions of Theorem 2.1. 

Conversely, if $R(z)$ is the corresponding multiplier then $R(z)$ is a product of $n$ elementary 
multipliers with one $\kappa$ equal to $-1$ and one $\delta$ equal to $+1$. The explicit construction of the proof of Theorem 2.1 shows that all these multipliers are polynomials, hence, $R(z)$ is a polynomial. The fact that $\delta_1=\dots=\delta_m$ implies that $R(z)$ is a linear polynomial of the form $z-B$ for some $B\in\Mat$ (to see that, it suffices to look at the asymptotics of the canonical solutions). 
We have 
$$
A(z)=R^{-1}(z+1)\wt A(z) R(z)=(z+I-B)^{-1}\wt A(z) (z-B).
$$
Comparing the determinants of both side we conclude that $Sp(B)=\{a^{(i)}_1,\dots,a^{(i)}_m\}$. Since no two eigenvalues are different by an integer, $B$ and $B-I$ have no common eigenvalues. This implies that $(z+I-B)^{-1}\wt A(z)$ must be a polynomial, and hence $z-B$ is a right divisor of $A(z)$.\qed
\enddemo

For any $k=(k_1,\dots,k_n)\in \Z^n$ we introduce matrices $B_1(k),\dots,B_n(k)$ such that the right divisors of $\s(k_1,\dots,k_n)A(z)$ have the form
$
z-B_i(k)
$
with 
$$
Sp(B_i(k))=\{a^{(i)}_1-k_i,\dots,a^{(i)}_n-k_i\},\quad i=1,\dots,n.
$$
They are defined for generic $A(z)$ from the varieties $\M(\cdots)$ introduced in the previous section.

\proclaim{Proposition 3.7 (difference Schlesinger equations)} The matrices $\{B_i(k)\}$ (whenever exist) satisfy the following equations:
$$
\gather
B_i(\dots)-B_i(\dots ,k_j+1,\dots)=
B_j(\dots)-B_j(\dots ,k_i+1,\dots),
\tag 3.1\\
{B_j(\dots ,k_i+1,\dots)B_i(\dots)=B_i(\dots ,k_j+1,\dots)B_j(\dots)},
\tag 3.2\\
B_i(k_1+1,\dots,k_n+1)=A_0^{-1}B_i(k_1,\dots,k_n)A_0-I,
\tag 3.3
\endgather 
$$ 
where $i,j=1,\dots,n$, and dots in the arguments mean that other $k_l$'s remain unchanged.
\endproclaim
\example{Remark 3.8} The first two equations above are equivalent to 
$$
\gathered
\bigl(z-B_i(\dots,k_j+1,\dots)\bigr)\bigl(z-B_j(\dots)\bigr)=
\bigl(z-B_j(\dots,k_i+1,\dots)\bigr)\bigl(z-B_i(\dots)\bigr).
\endgathered
\tag 3.4
$$
\endexample
\demo{Proof of Proposition 3.7} The uniqueness part of Theorem 2.1 implies that
$$
\gathered
\s(0,\dots,0,\overset{(i)}\to 1,0,\dots,0)\circ\s(0,\dots,0,\overset{(j)}\to 1,0,\dots,0)\circ \s(k_1,\dots,k_n)\\=
\s(0,\dots,0,\overset{(j)}\to 1,0,\dots,0)\circ\s(0,\dots,0,\overset{(i)}\to 1,0,\dots,0)\circ \s(k_1,\dots,k_n).
\endgathered
$$
Thus, the corresponding products of the multipliers are equal, which gives \tht{3.4}. This proves \tht{3.1}, \tht{3.2}. The relation \tht{3.3} follows from the fact that the multiplier for $\s(1,\dots,1)$ is equal 
to $A_0^{-1}A(z)$, and $\wt A(z)=\s(1,\dots,1)=A_0^{-1}A(z+1)A_0$. This means that the right divisors for $\wt A(z)$ can be obtained from those for $A(z)$ by shifting $z$ by 1 and conjugating by $A_0$.\qed
\enddemo

\proclaim{Theorem 3.9} Fix $mn$ complex numbers $\bigl\{a_j^{(i)}\bigr\}_{i=1,j=1}^{n,m}$ such that no two of them are different by an integer, and an integer $M>0$. Then for generic $B_1,\dots,B_n\in\Mat$, $Sp(B_i)=\bigl\{a_j^{(i)}\bigr\}_{j=1}^m$, there exists a unique solution 
$$
\{B_i(k_1,\dots,k_n): \max_{i=1,\dots,n}|k_i|\le M\}
$$ 
of the difference Schlesinger equations \tht{3.1}-\tht{3.3} with 
$$
Sp(A_0)=\{\rho_1,\dots,\rho_n\},\quad \rho_i/\rho_j\notin\Bbb R\,\,\text{  for  }
\,\,i\ne j,\quad \rho_i\ne 0\,\,\text{  for  }\,\, i=1,\dots,n,
$$
such that 
$$
Sp(B_i(k_1,\dots,k_n))=Sp(B_i)-k_i \quad\text{and}\quad B_i(0,\dots,0)=B_i \quad\text{for all  }\,\, i=1,\dots,n.
$$
The matrix elements of $B_i(k)$ are rational functions of the matrix elements of the initial conditions $\{B_i\}_{i=1}^n$. Moreover, these rational functions do not depend on the eigenvalues $\bigl\{a_i^{(j)}\bigr\}$.
\endproclaim

\example{Remark 3.10} As we will see later, this theorem also extends to the case of arbitrary invertible $A_0$.
\endexample

\demo{Proof} The existence and rationality of the flows have already been proved. Indeed, without loss of generality we can assume that $A_0$ is diagonal (the equations \tht{3.1}-\tht{3.3} remain intact if we conjugate all $B_i(k)$ and $A_0$ by the same constant matrix). By Lemma 3.4 we can construct a (unique) degree $n$ polynomial $A(z)$ with the highest coefficient $A_0$, such that $\{B_i\}$ is the set of its right divisors.
Then, using Theorem 2.1, we can define $\s(k)$ and hence $\{B_i(k)\}$. By Proposition 3.7 they will satisfy \tht{3.1}-\tht{3.3}. Moreover, all operations involved in this construction are rational. 

Thus, it remains to prove uniqueness and the fact that the rational functions involved do not depend on the eigenvalues. A simple computation shows that for any $X,Y,S,T\in\Mat$, the relation $(z-X)(z-Y)=(z-S)(z-T)$ implies 
$$
\gather
Y=(X-S)^{-1}S(X-S),\qquad T=(X-S)^{-1}X(X-S),\tag 3.5\\
X=(Y-T)T(Y-T)^{-1},\qquad S=(Y-T)Y(Y-T)^{-1},\tag 3.6
\endgather
$$
whenever the corresponding matrices are invertible, cf. \cite{GRW}.
Applying this observation to \tht{3.4}, we see that, generically, $\{B_i=B_i(0,\dots,0)\}$ uniquely define all 
$$
B_i(\epsilon^{(i)}_1,\dots,\epsilon^{(i)}_{i-1},\overset {(i)}\to0,\epsilon^{(i)}_{i+1},\dots,\epsilon^{(i)}_n),\qquad \epsilon^{(i)}_j=0,1.
\tag 3.7
$$
Moreover, they are all given by rational expressions involving the initial conditions $\{B_i\}$ only. To move further, we need the following lemma.
\proclaim{Lemma 3.11} For generic $X,Y\in \Mat$ with fixed disjoint spectra, there exist unique $S,T\in\Mat$ such that 
$$
(z-X)(z-Y)=(z-S)(z-T),\qquad
Sp(S)=Sp(Y),\quad Sp(T)=Sp(X).
$$
The matrix elements of $S$ and $T$ are rational functions of the matrix elements of $X$ and $Y$ which do not depend on the spectra of $X$ and $Y$.
\endproclaim
\demo{Proof 1} Lemma 3.3 proves the uniqueness and shows how to construct $T$ is we know the eigenvalues $x_1,\dots,x_m$ of $X$ and vectors $v_i$ 
such that $(x_i-X)(x_i-Y)v_i=0$. If we normalize $v_i$'s in the same way, for example, by requiring the first coordinate to be equal to 1 (this can be done generically), then using the construction of Lemma 3.3 we obtain the matrix elements of $T$ as rational functions in the matrix elements of $X,Y$ and
$x_1,\dots,x_n$. However, it is easy to see that these rational functions are symmetric with respect to the permutations of $x_1,\dots,x_n$, which means that they depend only on the elementary symmetric functions $\sum_{i_1<\dots<i_k}x_{i_1}\cdots x_{i_k}$ of $x_i$'s. But these are the coefficients of the characteristic polynomial of $X$, and hence they are expressible as polynomials in the matrix elements of $X$. \qed
\enddemo
\demo{Proof 2, see \cite{O}} The uniqueness follows from Lemma 3.3. To prove the existence, denote by $\Lambda$ the solution of the equation $Y\Lambda-\Lambda X=I$. Generically, it exists, it is unique and invertible. Set 
$$
S=X+\Lambda^{-1},\qquad T=Y-\Lambda^{-1}.
$$
Then it is easy to see that $(z-X)(z-Y)=(z-S)(z-T)$. Furthermore, if $Y$ and $T$ have a common eigenvalue then they must have a common eigenvector, which contradicts the invertibility of $Y-T=\Lambda^{-1}$. Hence, $Sp(T)=Sp(X)$ and $Sp(S)=Sp(Y)$.\qed
\enddemo

\example{Remark 3.12} In the case of $2$ by $2$ matrices, it is not hard to produce an explicit formula for $S$ and $T$ in terms of $X$ and $Y$:
$$
S=(X+Y-\operatorname{Tr}Y)Y(X+Y-\operatorname{Tr}Y)^{-1},\quad
T=(X+Y-\operatorname{Tr}X)^{-1}X(X+Y-\operatorname{Tr}X).
\tag 3.8
$$
\endexample

Now let us return to the proof of Theorem 3.9. Recall that we already proved that the initial conditions define \tht{3.7} uniquely. Now let us use \tht{3.4} with 
$$
(k_1,\dots,k_n)=(1,\dots,1,\overset{(j)}\to 0,1,\dots,1),\quad j\ne i.
$$
By \tht{3.3}, we know what $B_i(1,\dots,1)$ is. Thus, we know both matrices on the left--hand side of \tht{3.4}, and hence, by Lemma 3.11, we can compute both matrices on the right--hand side of \tht{3.4}, in particular, $B_i(1,\dots,1,\overset{(j)}\to 0,1,\dots,1)$.

Now take \tht{3.4} with 
$$
(k_1,\dots,k_n)=(1,\dots,1,\overset{(j)}\to 0,1,\dots,1,\overset{(l)}\to 0,1,\dots,1),
$$
where $i,j,l$ are pairwise distinct.
Applying Lemma 3.11 again, we find all 
$$
B_i(1,\dots,1,\overset{(j)}\to 0,1,\dots,1,\overset{(l)}\to 0,1,\dots,1).
$$
Continuing the computations in this fashion (changing one more 1 to 0 in $(k_1,\dots,k_n)$ on each step), we obtain all
$$
B_i(\epsilon^{(i)}_1,\dots,\epsilon^{(i)}_{i-1},\overset {(i)}\to 1,\epsilon^{(i)}_{i+1},\dots,\epsilon^{(i)}_n),\qquad \epsilon^{(i)}_j=0,1.
$$
Together with \tht{3.3} (and \tht{3.7}) this computes all $B_i(k)$ with
$\max|k_i|\le 1$. Iterating this procedure, we complete the proof.\qed
\enddemo

\head 4. An alternative description of the Schlesinger flows.
\endhead

The goal of this section is to provide yet another set of coordinates for the polynomials $A(z)$, in which the flows described in the previous section can be easily defined. In particular, this will lead to a different proof of Theorem 3.9, which will be valid for an arbitrary invertible $A_0$.

\proclaim{Proposition 4.1} In the assumptions of Theorem 3.9, the monic degree $n$ polynomial
$$
(z-B_1(0,1,\dots,1))(z-B_2(0,0,1,\dots,1))\cdots(z-B_n(0,\dots,0))
$$
has $z-B_i$, $i=1,\dots,n$, as its right divisors.
\endproclaim
This statement and Theorem 3.9 provide a proof for Remark 3.5(1).
\demo{Proof} Using \tht{3.4} we obtain $(j>i)$
$$
\gathered
\left(z-B_i(0,\dots,0,\overset{(j)}\to 1,\dots,1)\right)\left(z-B_j(0,\dots,0,\overset{(j+1)}\to 1,\dots,1)\right)\\=
\left(z-B_j(0,\dots,0,\overset{(i)}\to 1,0,\dots,0,\overset{(j+1)}\to 1,\dots,1)\right)\left(z-B_i(0,\dots,0,\overset{(j+1)}\to 1,\dots,1)\right).
\endgathered
$$
Using this commutation relation, we can move the factor $(z-B_i(\cdots))$
in the product above, to the right most position, where it will turn into $(z-B_i(0,\dots,0))=(z-B_i)$.\qed
\enddemo

Let us introduce the notation ($l_1,\dots,l_n\in\Z$)
$$
C_i(l_1,\dots,l_n):=B_i(l_1,\dots,l_i,l_{i+1}+1,\dots,l_n+1),\qquad i=1,\dots,n.
$$
If we denote by $A(z)$ the polynomial of degree $n$ with highest coefficient $A_0$ such that $\{z-B_i\}$ are its left divisors, then
the definition of $B_i(k)$ and Proposition 4.1 imply that for any
$l=(l_1,\dots,l_n)\in\Z^n$,
$$
\s(l_1,\dots,l_n)A(z)=A_0\bigl(z-C_1(l)\bigr)\cdots
\bigl(z-C_n(l)\bigr).
\tag 4.1
$$
(To apply Proposition 4.1, we also used an easy fact that for any solution $\{B_i(k)\}$ of \tht{3.1}-\tht{3.3} and any $l_1,\dots,l_n\in\Z$,
$$
B_i'(k_1,\dots,k_n):=B_i(k_1+l_1,\dots,k_n+l_n),\qquad i=1,\dots,n,
$$
also form a solution of \tht{3.1}-\tht{3.3}.)

\proclaim{Lemma 4.2} The map $\{B_i\}\mapsto\{C_i\}$ is birational. 
\endproclaim
\demo{Proof} The rationality of the forward map follows from Theorem 3.9. The rationality of the inverse map follows from Lemma 3.3 (indeed, we just need to find the right divisors of the known matrix $\s(l_1,\dots,l_n)A(z)$).
Even though it looks like to construct $B_i$ we need to know the eigenvalues of $C_i$, it is clear that by normalizing the eigenvectors of $A(z)$ corresponding to these eigenvalues, in the same way, we will obtain a formula for $B_i$ which will be symmetric with respect to the permutations of these eigenvalues, and thus we can rewrite it through the matrix elements of $C_i$'s only (this argument was already used in the first proof of Lemma 3.11 above). 
\qed
\enddemo

Our goal is to describe the transformations $\s(k)$ in terms of $\{C_i\}$. We need a preliminary lemma which generalizes Lemma 3.11.

\proclaim{Lemma 4.3} For generic $X_1,\dots,X_N\in\Mat$ with fixed disjoint spectra and any permutation $\sigma\in S_N$, there exist unique $Y_1,\dots,Y_N\in \Mat$ such that
$Sp(Y_i)=Sp(X_i)$ for all $i=1,\dots,N$, and
$$
(z-X_1)\cdots(z-X_N)=(z-Y_{\sigma(1)})\cdots(z-Y_{\sigma(N)}).
$$
The matrix elements of $\{Y_i\}$ are rational functions of the matrix elements of $\{X_i\}$ which do not depend on the spectra of $\{X_i\}$.
\endproclaim
\demo{Proof} The existence and rationality claims follow from Lemma 3.11, because elementary transpositions $(i,i+1)$ generate the symmetric group $S_N$. To show uniqueness, we rewrite the equality
$$
(z-Y_1')\cdots(z-Y_N')=(z-Y_1'')\cdots(z-Y_N''),\qquad Sp(Y_i')=Sp(Y_i''),
$$
in the form
$$
(z-Y_1'')^{-1}(z-Y_1')=\Bigl((z-Y_2'')\cdots(z-Y_N'')\Bigr)
\Bigl((z-Y_2')\cdots(z-Y_N')\Bigr)^{-1}.
$$
If the spectrum of $Y_1''$ is disjoint with the spectra of $Y_2',\dots,Y_m'$, then both sides of the last equality are entire because they cannot possibly have common poles. Since both sides tend to $I$ as $z\to\infty$, by Liouville's theorem we conclude that both sides are identically equal to $I$, and  $Y_1'=Y_1''$. Induction on $N$ concludes the proof. \qed 
\enddemo

\proclaim{Proposition 4.4} In the assumptions of Theorem 3.9, $\{C_i(l)\}$ satisfy the equations
$$
\gathered
\bigl(z+1-C_i\bigr)
\cdots
\bigl(z+1-C_n\bigr)
A_0\bigl(z-C_1\bigr)\cdots
\bigl(z-C_{i-1}\bigr)
\\=
\bigl(z+1-\wt C_{i+1}\bigr)
\cdots
\bigl(z+1-\wt C_n\bigr) 
A_0\bigl(z-\wt C_1\bigr)\cdots
\bigl(z-\wt C_{i}\bigr),
\\
C_j=C_j(l_1,\dots,l_n),\quad \wt C_j=C_j(l_1,\dots,l_{i-1},l_i+1,l_{i+1},\dots,l_n)\text{  for all  }j,
\endgathered
\tag 4.2
$$
and
$$
C_i(l_1+1,\dots,l_n+1)=A_0^{-1}C_i(l_1,\dots,l_n)A_0-I.
\tag 4.3
$$
In both equations $i=1,\dots,n$ is arbitrary.
\endproclaim
\demo{Proof} The relation \tht{4.3} is a direct corollary of \tht{3.3}. Let us prove \tht{4.2}. Proposition 3.6 implies that the multiplier for the shift 
$$
(l_1,\dots,l_n)\mapsto (l_1,\dots,l_{i-1},l_i+1,l_{i+1},\dots,l_n)
$$
has the form $R(z)=z-B_i(l)$. Thus, \tht{4.1} gives
$$
A_0(z-\wt C_1)\cdots (z-\wt C_n)(z-B_i(l))=
(z+1-B_i(l))A_0(z-C_1)\cdots (z-C_n) .
$$
Comparing the spectra of factors on both sides and applying Lemma 4.3, we get
$$
\gathered
A_0(z-\wt C_1)\cdots (z-\wt C_i)=(z+1-B_i(l))A_0(z-C_1)\cdots(z-C_{i-1}),\\
(z-\wt C_{i+1})\cdots (z-\wt C_n)(z-B_i(l))=
(z-C_i)\cdots (z-C_n).
\endgathered
$$
Combining these two relations and shifting $z\mapsto z+1$ in the second one, we arrive at \tht{4.2}.\qed 
\enddemo
\proclaim{Theorem 4.5} Fix $mn$ complex numbers $\bigl\{a_j^{(i)}\bigr\}_{i=1,j=1}^{n,m}$ such that no two of them are different by an integer, an integer $M>0$, and any nondegenerate $A_0\in\Mat$. Then for generic $C_1,\dots,C_n\in\Mat$, $Sp(C_i)=\bigl\{a_j^{(i)}\bigr\}_{j=1}^m$, there exists a unique solution
$$
\{C_i(l_1,\dots,l_n): \max_{i=1,\dots,n}|l_i|\le M\}
$$ 
of the equations \tht{4.2} and, consequently, \tht{4.3}, such that
$$
Sp(C_i(l_1,\dots,l_n))=Sp(C_i)-l_i \quad\text{and}\quad C_i(0,\dots,0)=C_i \quad\text{for all  }\,\, i=1,\dots,n.
$$
The matrix elements of $C_i(l)$ are rational functions of the matrix elements of the initial conditions $\{C_i\}_{i=1}^n$, and these rational functions do not depend on the eigenvalues $\bigl\{a_i^{(j)}\bigr\}$.
Moreover, 
$$
B_i(k_1,\dots,k_n):=C_i(k_1,\dots,k_i,k_{i+1}-1,\dots,k_n-1),\qquad i=1,\dots,n,
$$
solve the difference Schlesinger equations \tht{3.1}-\tht{3.3}.
\endproclaim
\example{Remark 4.6} If the ratios of eigenvalues of $A_0$ are not real then Theorem 3.9, Lemma 4.2, and Proposition 4.4 provide a proof of Theorem 4.5. However, our goal is to provide an independent proof of this theorem, thus giving a different proof of Theorem 3.9 with arbitrary invertible $A_0$, cf. Remark 3.10.
\endexample

To prove Theorem 4.5 we will develop a rather general formalism.

\subhead (a) Semigroup\endsubhead
Let $P$ be a semigroup and $P_0$ be its subset. We assume that every element of $P_0$ has a {\it type}. The types of two different elements $p_1,p_2\in P_0$ may be the same, which will denoted by $t(p_1)=t(p_2)$, and may be {\it disjoint}, which will be denoted by $t(p_1)\perp t(p_2)$. The types may also be neither equal nor disjoint.  

\proclaim{Assumption 4.7} For any elements $p_1,\dots,p_N\in P_0$ such that their types are pairwise disjoint:
$$
t(p_i)\perp t(p_j),\qquad i\ne j,\quad 1\le i,j\le N,
$$
and for any permutation $\sigma\in S_N$ there exist unique elements
$\hat p_1,\dots \hat p_N\in P_0$ such that $t(\hat p_i)=t(p_i)$, $i=1,\dots,N$, and 
$$
p_1\cdots p_N=\hat p_{\sigma(1)}\cdots \hat p_{\sigma(N)}.
$$
\endproclaim

We will be interested in the situation when 
$$
\gathered
P=P^{\Mat}=\left\{z^k+Q_{1}z^{k-1}+\dots+Q_k\,|\, Q_i\in\Mat,
\,1\le i\le k\right\},\\
P_0=P_0^{\Mat}=\left\{z-Q\,|\,Q\in\Mat\right\},
\\
t(z-Q)=\{z\in\C\,|\,\det(z-Q)=0\}=Sp(Q).
\endgathered
$$

The notions of equality and disjointness for types are the natural ones for the $m$-point subsets of $\C$. Lemma 4.3 shows that $P^{\Mat}$ and $P_0^{\Mat}$ satisfy Assumption 4.7 generically. 

\proclaim{Proposition 4.8} Let $P$ be a semigroup satisfying the Assumption 4.7. Assume that we have an equality in $P$ of the form
$$
\left(p_1^{(1)}\cdots p_{m_1}^{(1)}\right)
\cdots \left(p_1^{(k)}\cdots p_{m_k}^{(k)}\right)=
\left(q_1^{(1)}\cdots q_{m_1}^{(1)}\right)
\cdots \left(q_1^{(k)}\cdots q_{m_k}^{(k)}\right),
\tag 4.4
$$
where all $p_i^{(j)},q_i^{(j)}$ are from $P_0$, the types of all elements on the left hand side are pairwise disjoint, the types of all elements on the right-hand side are pairwise disjoint, and 
$$
\left\{t(p_1^{(j)}),\dots,t(p_{m_j}^{(j)})\right\}=
\left\{t(q_1^{(j)}),\dots,t(q_{m_j}^{(j)})\right\}
$$
for all $j=1,\dots,k$. Then
$$
p_1^{(j)}\cdots p_{m_j}^{(j)}=
q_1^{(j)}\cdots q_{m_j}^{(j)},\qquad j=1,\dots,k.
$$
\endproclaim
\demo{Proof} By Assumption 4.7, for any $j=1,\dots,k$, there exist 
$\hat q_1^{(j)},\dots, \hat q_{m_j}^{(j)}$ such that
$$
q_1^{(j)}\cdots q_{m_j}^{(j)}=\hat q_1^{(j)}\cdots \hat q_{m_j}^{(j)}
$$
and $t(p_i^{(j)})=t(\hat q_i^{(j)})$. Then by the uniqueness part of Assumption 4.7 applied to \tht{4.4} we obtain $p_i^{(j)}=\hat q_i^{(j)}$ for all $i,j$. This immediately implies the claim.\qed
\enddemo
We, essentially, used Proposition 4.8 in the proof of Proposition 4.4 above.
\subhead (b) Commuting flows on sequences\endsubhead
Denote by $\Cal P$ the set of all sequences $\{p_k\}_{k\in\Z}\subset P_0$ such that the types of all elements of a sequence are pairwise disjoint.

Fix an integer integer $n>0$. For any $l\in\Z$ we define a map $F_l:\Cal P\to\Cal P$ as follows:
$$
\gathered 
F_l:\{p_k\}_{k\in\Z}\mapsto \{q_k\}_{k\in \Z},\\
p_{l+\mu n}p_{l+\mu n+1}\cdots p_{l+(\mu+1)n-1}=
q_{l+\mu n+1}q_{l+\mu n+2}\cdots q_{l+(\mu+1)n},\quad \mu\in\Z,\\
t(q_j)=\cases t(p_{j-n}),&\text{if  }(j-l) \text{  divides  } n,\\
              t(p_{j}),& \text{otherwise}.
	\endcases
\endgathered
$$	          
In this definition $\mu$ ranges over all integers, and for each $\mu$ we use Assumption 4.7 for $\sigma = (12\cdots n)\in S_n$. Clearly, $F_l$ is invertible.

It is convenient to denote
$$
F_{l_1}\Bigl(F_{l_2}\bigl(\cdots F_{l_m}\left(\{p_k\}\right)\bigr)\Bigr)
=\{p_k^{l_1,\dots,l_m}\}.
$$
Then the second line in the definition above takes the form
$$
p_{l+\mu n}p_{l+\mu n+1}\cdots p_{l+(\mu+1)n-1}=
p^l_{l+\mu n+1}p^l_{l+\mu n+2}\cdots p^l_{l+(\mu+1)n}.
$$

For example, for $n=2$ we have
$$
\gathered
p_{2s-1}p_{2s}=p^1_{2s}p^1_{2s+1}, \qquad 
t(p^1_{2s})=t(p_{2s}),\quad t(p^1_{2s+1})=t(p_{2s-1}),\quad s\in\Z,\\
p_{2s}p_{2s+1}=p^2_{2s+1}p^2_{2s+2}, \qquad 
t(p^2_{2s+1})=t(p_{2s+1}),\quad t(p^2_{2s+2})=t(p_{2s}),\quad s\in\Z.
\endgathered
$$

It is immediately seen from the definition that $F_{l+\mu n}=F_l$ and
$$
p_k^{l,l+1,\dots,l+n-1}=p_{k-n}
\tag 4.5
$$
for any $k,l\in\Z$.
\proclaim{Theorem 4.9} {\rm (i)} For any $i,j\in\Z$, $F_i$ and $F_j$ commute. That is 
$$
p_k^{i,j}\equiv p_k^{j,i}\text{  for any  }\{p_k\}\in\Cal P.
$$

{\rm (ii)} For any $\{q_k\}\in\Cal P$ and any $i,j\in\Z$ such that
$0<j-i<n$, set $p_k=q_k^{i+1,\dots,j}$. Then 
$$
q_i^j\,p_j=p_j^i\,q_i.
\tag 4.6
$$
\endproclaim

\example{Remark 4.10} Part (i) of this theorem means that we have defined an action of $\Z^n$ on $\Cal P$. There is a much larger group which acts on $\Cal P$. Let $\pi:\Z\to\Z$ be a bijection such that for any $k\in\Z$ the sets
$$
I_k=\{i\in\Z: i<k,\,\pi(i)>\pi(k)\},\qquad
J_k=\{j\in\Z: j>k,\,\pi(j)<\pi(k)\}
$$
are finite: $I=\{i_1,\dots,i_s\}$, $J=\{j_1,\dots,j_t\}$. Then, given a sequence $\{p_k\}\in\Cal P$, we define $\{p^\pi_k\}\in\Cal P$ by
$$
p_{i_1}\cdots p_{i_s}\,\, p_k\,\, p_{j_1}\cdots p_{j_t}=
p'_{j_1}\cdots p'_{j_t}\,\, p_{\pi(k)}^\pi\,\,p'_{i_1}\cdots p'_{i_s}
$$
where $t(p_l')=t(p_l)$ and $t\bigl(p_{\pi(k)}^\pi\bigr)=t(p_k)$. One can show that this defines an action of the group of all $\pi$ satisfying the condition above on the space $\Cal P$.  
The maps $\{F_l\}$ correspond to shifts by $n$ along $n$ nonintersecting arithmetic progressions $\{l+\mu n:\mu\in\Z\}$, hence they must commute.
\endexample

\demo{Proof of Theorem 4.9}
Since $F_l=F_{l+n}$, it suffices to assume that $0<j-i<n$. Consider the product
$$
\Pi= p_ip_{i+1}\dots p_{j+2n-1}.
$$
On one hand, we have
$$
\gathered
\Pi=\bigl(p_i\cdots p_{i+n-1}\bigr)\bigl(p_{i+n}\cdots p_{i+2n-1}\bigr)
p_{i+2n}\cdots p_{j+2n-1}\\
=\bigl(p_{i+1}^i\cdots p_{i+n}^i\bigr)
\bigl(p_{i+n+1}^i\cdots p_{i+2n}^i\bigr)
p_{i+2n}\cdots p_{j+2n-1}\\
=p_{i+1}^i\cdots p_{j-1}^i\bigl(p_j^i\cdots p_{j+n-1}^i)
p_{j+n}^i\cdots p_{i+2n}^i\,p_{i+2n}\cdots p_{j+2n-1}\\
=p_{i+1}^i\cdots p_{j-1}^i\,p_{j+1}^{j,i}\cdots p_{j+n}^{j,i}\,
p_{j+n}^i\cdots p_{i+2n}^i\,p_{i+2n}\cdots p_{j+2n-1}.
\endgathered
$$
On the other hand, we have
$$
\gathered
\Pi=p_i\cdots p_{j-1}\bigl(p_j\cdots p_{j+n-1}\bigr)
\bigl(p_{j+n}\cdots p_{j+2n-1}\bigr)\\
=p_i\cdots p_{j-1}\bigl(p_{j+1}^j\cdots p_{j+n}^j\bigr)
\bigl(p_{j+n+1}^j\cdots p_{j+2n}^j\bigr)\\
=p_i\cdots p_{j-1}\,p_{j+1}^j\cdots p_{i+n-1}^j\bigl(p_{i+n}^j\cdots 
p_{i+2n-1}^j\bigr)p_{i+2n}^j\cdots p_{j+2n}^j\\
=p_i\cdots p_{j-1}\,p_{j+1}^j\cdots p_{i+n-1}^j\,p_{i+n+1}^{i,j}\cdots 
p_{i+2n}^{i,j}\,p_{i+2n}^j\cdots p_{j+2n}^j.
\endgathered
$$
Thus, we obtain
$$
\gathered
\bigl(p_{i+1}^i\cdots p_{j-1}^i\,p_{j+1}^{j,i}\cdots p_{i+n}^{j,i}\bigr)\bigl(p_{i+n+1}^{j,i}\cdots p_{j+n}^{j,i}\bigr)
\bigl(
p_{j+n}^i\cdots p_{i+2n}^i\,p_{i+2n}\cdots p_{j+2n-1}\bigr)\\
=\bigl(p_i\cdots p_{j-1}\,p_{j+1}^j\cdots p_{i+n-1}^j\bigr)\bigl( p_{i+n+1}^{i,j}\cdots p_{j+n}^{i,j}\bigr)\bigl(p_{j+n+1}^{i,j}\cdots 
p_{i+2n}^{i,j}\,p_{i+2n}^j\cdots p_{j+2n}^j\bigr).
\endgathered
\tag 4.7
$$
Comparing the types in the three factors on the left and on the right, we see that we are in a position to apply Proposition 4.8. It implies, in particular, that the middle factors are equal. Since the order of the types in the middle factors is the same, these middle factors must be equal termwise:
$$
p_k^{j,i}=p_k^{i,j},\qquad i+n+1\le k\le j+n.
$$
Since $F_l=F_{l+n}$ for all $n$, and $i$ and $j$ are arbitrary, we see that $p_k^{j,i}=p_k^{i,j},\quad i+1 \le k\le j.$
Switching from $(i,j)$ to $(j,i+n)$, we get $p_k^{i,j}=p_k^{j,i}$ for 
$j+1\le k\le i+n$. Thus, the commutativity relation is proved for $i+1\le k\le i+n$, and thus for all $k\in\Z$. The proof of the first part of Theorem 4.9 is complete.

In order to prove Theorem 4.9(ii), we need to compare the first and the third factors of the two sides of \tht{4.7}. The first factors give
$$
p_{i+1}^i\cdots p_{j-1}^i\,p_{j+1}^{j,i}\cdots p_{i+n}^{j,i}=p_i\cdots p_{j-1}\,p_{j+1}^j\cdots p_{i+n-1}^j.
$$
Commuting $p_i$ in the right--hand side to the right and using Proposition 4.8, we see that
$$
p_i\cdots p_s=p_{i+1}^i\cdots p_s^i\, x_{s,i}
\tag 4.8
$$
where $i+1\le s\le j-1$, $x_{s,i}\in P_0$, and $t(x_{s,i})=t(p_i)$. Since $j$ is arbitrary (but $j-i<n$), we can assume that $i+1\le s\le i+n-2$.
Note also that \tht{4.8} also holds for $s=i+n-1$ with $x_{i+n-1,i}=p_{i+n}^i$, as follows from the definition of $F_i$.

Similarly, looking at the third factors and substituting $\{p_k\}$ by $F_j^{-1}\{p_k\}$, we get
$$
p_t\cdots p_j=y_{t,j}\, \left(F_j^{-1}\{p_k\}\right)_t\cdots \left(F_j^{-1}\{p_k\}\right)_{j-1}
\tag 4.9
$$
where $j-n+2\le t\le j-1$, $y_{t,j}\in P_0$, and $t(y_{t,j})=t(p_j)$.
Again, this also holds for $t=j-n+1$ with $y_{j+n-1,j}=\left(F_j^{-1}\{p_k\}\right)_{j-n}$.

\proclaim{Lemma 4.11} For $i+1\le s\le i+n-1$, we have 
$$
x_{s,i}=\left(F^{-1}_{i+1}\circ\cdots\circ F_s^{-1}\{p_k\}\right)_i.
$$
\endproclaim
\demo{Proof} Induction on $s-i$.  
To prove both the base $s=i+1$ of the induction and the induction step we first use \tht{4.9} to write
$$
p_i\cdots p_s=y_{i,s}\, \left(F_s^{-1}\{p_k\}\right)_i\cdots \left(F_s^{-1}\{p_k\}\right)_{s-1}
$$
and now use the induction hypothesis on the factors which go after $y_{i,s}$ to obtain
$$
p_i\cdots p_s=y_{i,s}\,
\left(F_iF_s^{-1}\{p_k\}\right)_{i+1}\cdots 
\left(F_iF_s^{-1}\{p_k\}\right)_{s-1}
\left(F_{i+1}^{-1}\circ\cdots \circ F_{s}^{-1}\{p_k\}\right)_i.
$$
(If $s=i+1$ then the second step is empty.)
Since 
$$
t\left(\left(F_{i+1}^{-1}\circ\cdots\circ F_{s}^{-1}\{p_k\}\right)_i\right)=t(x_{s,i})=t(p_i),
$$
comparing with \tht{4.8} we conclude that $x_{s,i}=\left(F^{-1}_{i+1}\circ\cdots\circ F_s^{-1}\{p_k\}\right)_i$. 

This argument works for $s\le i+n-2$. For $s=i-n+1$ the lemma follows from \tht{4.5}. \qed
\enddemo

Now we return to the second part of Theorem 4.9. Applying Lemma 4.11 to all but one factors in $p_i\cdots p_j$, and then to all factors in $p_i\cdots p_j$, we obtain
$$
\gathered
p_i\cdots p_j=p_{i+1}^i\cdots p_{j-1}^i\left(F^{-1}_{i+1}\circ\cdots\circ F_{j-1}^{-1}\{p_k\}\right)_ip_j\\=
p_{i+1}^i\dots p_{j-1}^ip_j^i\left(F_{i+1}^{-1}\circ\cdots\circ F_{j}^{-1}\{p_k\}\right)_i.
\endgathered
$$
In the last two products all but last two factors coincide. By Proposition 4.8, this means that the products of last two ones also coincide:
$$
\left(F_{i+1}^{-1}\circ\cdots\circ F_{j-1}^{-1}\{p_k\}\right)_ip_j=
p_j^i\left(F_{i+1}^{-1}\circ\cdots\circ F_{j}^{-1}\{p_k\}\right)_i.
$$
Renaming $F_{i+1}^{-1}\circ\cdots\circ F_{j}^{-1}\{p_k\}$ by  $\{q_k\}$ we arrive at \tht{4.6}. The proof of Theorem 4.9 is complete. \qed
\enddemo

\subhead (c) Proof of Theorem 4.5 \endsubhead Let us concentrate on the case $P=P^{\Mat}$, $P_0=P_0^{\Mat}$, see (a) above. Since Assumption 4.7 generically holds in this case (see Lemma 4.3), we will be acting like it always holds, keeping in mind that all the claims we prove hold only generically. Set $p_i=z-C_i$, $i=1,\dots,n$, where $C_i=C_i(0,\dots,0)$ are as in Theorem 4.5. More generally, define 
$$
p_{i+\mu n}=z-\mu-A_0^\mu C_i(0,\dots,0) A_0^{-\mu},\qquad i=1,\dots,n,\quad\mu\in\Z,
\tag 4.10
$$
where $A_0$ is an arbitrary invertible element of $\Mat$. The assumption
that no two numbers of the set $\{a_j^{(i)}\}$ are different by an integer, guarantees that $\{p_k\}\in\Cal P$. Now define $\{C_i(l_1,\dots,l_n)\}$ by
$$
{\left(F_1^{l_1}\cdots F_n^{l_n}\{p_k\}\right)}_{i+\mu n}=
z-\mu-A_0^{\mu}C_i(l_1,\dots,l_n)A_0^{-\mu},\qquad i=1,\dots,n,\quad\mu\in\Z.
$$ 
(It is immediately seen that the subset of $\Cal P$ consisting of sequences $\{p_k=z-Q_k\}$ such that $Q_{k+n}=I+A_0Q_kA_0^{-1}$ is stable under the flows $F_1,\dots,F_n$, which shows that $C_i(l)$ are well-defined.) The very definition of $F_i$ implies \tht{4.2}. Furthermore, \tht{4.3} is a direct corollary of \tht{4.5}. It is easy that $Sp(C_i(l))=Sp(C_i)-l_i$, and this gives the existence part of Theorem 4.5. The uniqueness and rationality claims follow from Lemma 4.3.
Finally, let us show that
$$
B_i(k_1,\dots,k_n)=C_i(k_1,\dots,k_i,k_{i+1}-1,\dots,k_n-1)
$$
solve \tht{3.1}-\tht{3.3}. The relation \tht{3.3} is equivalent to \tht{4.3}. We will derive \tht{3.4} (and hence \tht{3.1}, \tht{3.2}) from Theorem 4.9(ii).

Fix $1\le i<j\le n$ and define 
$$
\{\tilde p_k\}=F_1^{k_1}\cdots F_n^{k_n}\{p_k\},\quad
\{\tilde q_k\}=F_{i+1}^{-1}\circ\dots\circ F_{j}^{-1}\{\tilde p_k\}.
$$
Then 
$$
\aligned
\tilde p_j&=z-C_j(k_1,\dots,k_n)=z-B_j(k_1,\dots,k_j,k_{j+1}+1,\dots,k_n+1),\\
\tilde 
q_i&=
z-B_i(k_1,\dots,k_{j},k_{j+1}+1,\dots,k_n+1),
\\
\tilde p_j^i&=z-B_j(k_1,\dots,k_{i-1},k_i+1,k_{i+1},\dots,k_{j},k_{j+1}+1,\dots,k_n+1),\\
\tilde q_i^j&=z-B_i(k_1,\dots,k_{j-1},k_{j}+1,\dots,k_n+1).
\endaligned
$$
If we apply the shift
$$
k_{j+1}\mapsto k_{j+1}-1,\dots ,k_n\mapsto k_n-1,
$$
then the equality $\tilde q_i^j\tilde p_j=\tilde p_j^i\tilde q_i$ turns into \tht{3.4}. This completes the proof of Theorem 4.5.\qed

\example{Remark 4.12}
The set of sequences $\{p_k=z-Q_k\}$ with $Q_{k+n}=I+A_0Q_kA_0^{-1}$ is also stable under the action of permutations $\pi:\Z\to\Z$ (see Remark 4.10) of the form
$$
\pi(i+\mu n)=\sigma(i)+\mu n, \qquad \sigma\in S_n,\quad i=1,\dots,n,\quad \mu\in\Z.
$$
Define $\hat C_i(l)=C_i(l)+l_iI$. Then we obtain a birational action of the semidirect product $\Z^n\ltimes S_n$ on 
$\{\hat C_1,\dots,\hat C_n\}\in (\Mat)^n$ which preserves the spectra of $\hat C_i's$.
\endexample

\example{Remark 4.13}
If instead of \tht{4.10} we use periodic initial conditions
$$
p_{i+\mu n}=z-C_i,\qquad i=1,\dots,n,\quad\mu\in\Z,
$$
which corresponds to the autonomous limit of the difference Schlesinger 
mentioned in the Introduction,
then the maps $F_1,\dots,F_n$ are exactly the {\it monodromy maps} constructed by Veselov \cite{V} in the framework of set-theoretical solutions of the quantum Yang-Baxter equation. We refer to \cite{V} for details and further references on the subject. 
\endexample

\example{Remark 4.14}
Solutions of the $q$-difference Schlesinger equations mentioned in the Introduction are obtained from considering $\{p_k=z-Q_k\}$ with $Q_{k+n}=qA_0Q_kA_0^{-1}$.
\endexample

\head 5. Continuous limit
\endhead

We start with a brief survey of the classical deformation theory for linear matrix differential equations, which is due to Riemann, Schlesinger, Fuchs, and Garnier, see \cite{JMU}, \cite{JM} for details.

Consider a first order matrix system of ordinary linear differential equations
$$
\frac{d\Y}{d\ze}=\B(\ze)\Y(\ze),\qquad \B(\ze)=\B_\infty+\sum_{k=1}^n\frac{\B_k}{\ze-x_k}\,.
\tag 5.1
$$

Here all matrices are in $\Mat$. We will assume that 
all $\B_k$'s can be diagonalized:
$$
\gathered
\B_k=G_kT_kG_k^{-1},\qquad T_k=\diag(t^{(k)}_1,\dots,t^{(k)}_m),\quad k=1,\dots,n,\\
\B_\infty=G_\infty\diag(s_1,\dots,s_n)G_\infty^{-1}
\endgathered
$$
where $t_i^{(k)}-t_j^{(k)}\notin \Z$, $i\ne j$, for all $k\ne\infty$, and $s_i\ne s_j$, $i\ne j$.

Alternatively, we may also consider
$$
\frac{d\Y}{d\ze}=\B(\ze)\Y(\ze),\qquad \B(\ze)=\sum_{k=1}^n\frac{\B_k}{\ze-x_k}\,,
\tag 5.2
$$
in which case we assume (in addition to the above assumption on $\B_k$, $k=1,\dots,m$) that
$$
-\sum_{k=1}^n\B_k=G_\infty T_\infty G_\infty^{-1},\qquad
T_\infty=\diag(t^{(\infty)}_1,\dots,t^{(\infty)}_m).
$$
with $t_i^{(\infty)}-t_j^{(\infty)}\notin \Z$ for $i\ne j$. 

Since we can conjugate $\Y$ and $\{\B_k\}$ by $G_\infty$, we may set $G_\infty=I$ without loss of generality.

One can show, see e.g. \cite{JMU, Proposition 2.1}, that
there exists a unique formal solution $\Y(\ze)$ of \tht{5.1} or \tht{5.2} of the form 
$$
\Y(\ze)=\hat\Y(\ze) \exp(T(\ze)),\qquad
\hat \Y(\ze)=I+\hat \Y_1\ze^{-1}+\hat \Y_2\ze^{-2}+\dots,
\tag 5.3
$$
where
$$
T(\ze)=\cases \diag(s_1,\dots,s_n) z+T_\infty\ln(z^{-1})&\quad
\text{for \tht{5.1}},\\T_\infty \ln(z^{-1})&\quad
\text{for \tht{5.2}},
\endcases
$$
with $T_\infty=\diag(t_1^{(\infty)},\dots,t_m^{(\infty)})$. (This formula is also the definition of $T_\infty$ for \tht{5.1}).
This is the analog of Propositions 1.1 and 1.2.

It turns out that for \tht{5.2} the series in \tht{5.3} is convergent, and after multiplication by $\exp(T(\ze))$ it defines a holomorphic near $\ze=\infty$ multi-valued function $\Y^\infty(\ze)$. However, in the case of \tht{5.1} this series is, generally speaking, divergent. Then the analog of Theorem 1.3 holds. Namely, there exist unique holomorphic solutions $\Y^{l,r}$ of \tht{5.1}, defined for $\Re\ze\ll 0$ and $\Re\ze \gg 0$, respectively, such that they have the asymptotic expansion \tht{5.3} as $\ze\to\infty$.\footnote{ As in the case of difference equations, one has to be careful in choosing the sector where $\ze$ may tend to $\infty$. One may always take $\arg\ze\in(\pi/2+\epsilon,3\pi/2-\epsilon)$ for $\Re \ze\ll 0$ and $\arg\ze\in(-\pi/2+\epsilon,\pi/2-\epsilon)$ for $\Re\ze\gg 0$. If $s_i-s_j\notin \Bbb R$, these sectors may be extended.} Since both these functions solve the same differential equation, there exist constant matrices $S^{\pm}$ such that the analytic continuations of $\Y^{l,r}$ in $\Im\ze\gg 0$ ($\ll 0$) are related by $\Y^l=\Y^r S^{\pm}$.  The matrices $S^{\pm}$ are called the {\it Stokes multipliers}.

It is also possible to determine the nature of solutions of \tht{5.1}, \tht{5.2} near the poles $x_1,\dots,x_k$. Namely, one can show that there exist locally holomorphic functions
$$
\hat \Y^{(k)}(\ze)=I+\Y^{(k)}_1(z-x_k)+\Y^{(k)}_2(z-x_k)^2+\dots
$$
such that for any solution $\Y(\ze)$, there exist constant matrices $C_k$ such that locally near $\ze=x_k$ we have
$$
\Y(\ze)=G_k\hat \Y^{(k)}(\ze)\exp\bigl(T_k\ln(z-x_k)\bigr)C_k.
$$
(Recall that $\{G_k\}$ were defined above by $\B_k=G_kT_kG_k^{-1}$.)
In particular, if we fix paths from $\ze=\infty$ (or $\pm\infty$ for \tht{5.1}), then we can define $\{C_k\}$ for the (analytic continuations of the) canonical solutions $\Y^\infty$ or $\Y^{l,r}$.

Thus, to any equation of the form \tht{5.1} or \tht{5.2}, we associate the following {\it monodromy data}: $\{T_k\}_{k=1}^n$ and $T_\infty$, $\{C_k\}_{k=1}^n$ computed 
for the canonical solution $\Y^\infty$ or $\Y^{l,r}$, and in the case of  \tht{5.1} we also add
the Stokes multipliers $S^{\pm}$ and the exponents $s_1,\dots,s_m$. 

If we analytically continue any solution $\Y(\ze)$ of \tht{5.1} or \tht{5.2} along a closed path $\gamma$ in $\C$ avoiding the singular points $\{x_k\}$ then the columns of $\Y$ will change into their linear combinations: $\Y\mapsto \Y M_\gamma$. Here $M_\gamma$ is a constant invertible matrix which depends only on the homotopy class of $\gamma$. It is called the {\it monodromy matrix} corresponding to $\gamma$. If $\gamma$ is a positive loop around $x_k$ then the corresponding monodromy matrix $M_k$ for the canonical solution $\Y^\infty$ or $\Y^{l,r}$ can be computed using the monodromy data introduced above:
$$
M_k=
C_k^{-1}\exp(2\pi i T_k)C_k.
$$

The basic problem of the {\it isomonodromy deformation} of the linear system \tht{5.1} or \tht{5.2} is to change $\B(\ze)$ in such a way that the monodromy data, or, more generally, the monodromy matrices $\{M_k\}$ remain invariant.

There are two types of isomonodromy deformations, both were discovered by Schlesinger \cite{Sch} and later generalized to the case of singularities of higher order in \cite{JMU}, \cite{JM}.

The first type is a continuous deformation which allows the singularities $x_1,\dots,x_n$ to move and describes $\{\B_k\}$ as functions of $x_j$'s.
This deformation leaves the whole set of monodromy data intact. The evolution of $\{\B_k\}$ is described by a system of partial differential equations called the {\it Schlesinger equations}:
$$
\frac{\partial \B_l}{\partial x_j}=\frac
{[\B_j,\B_l]}{x_j-x_l}\,,\qquad   \frac{\partial \B_j}{\partial
b_j}=\sum_{\Sb 1\le l\le n\\ l\ne j\endSb}\frac
{[\B_j,\B_l]}{x_l-x_j}-[\B_j,\B_\infty]\,, \qquad j,l=1,\dots,n,
\tag 5.4
$$
where for \tht{5.2} the term with $\B_\infty$ is dropped.

 It is not hard to show that this system has a local solution for arbitrary initial
conditions $\{\B_k(x_1^o,\dots,x_n^o)\}$. It is a much deeper fact (proved independently in \cite{Mal}, \cite{Miw}) that the Schlesinger equations with arbitrary initial conditions have a {\it global} meromorphic solution on the universal covering space of 
$$ 
\{(x_1,\dots,x_n)\in\C^n:x_i\ne x_j \text{  for  }i\ne j\}.
$$
To describe this fact, one often says that the system of Schlesinger equations enjoys the {\it Painlev\'e property}.

The second deformation (or, better to say, transformation) is an action of $\Z^{m(n+1)-1}$ on the space of $\B(z)$, which consists of multiplying $\Y(z)$ by an appropriate rational function on the left: $\Y(z)\mapsto \R(z)\Y(z)$. Such a transformation (called {\it Schlesinger transformation}) is uniquely determined by the shifts 
$$
t_j^{(k)}\mapsto t_j^{(k)}+\lambda_j^{(k)},\qquad k=1,\dots,n,\infty,
$$
of the eigenvalues of $\B_k$. Here all $\lambda_j^{(k)}$ are integers, and their total sum is equal to zero. Schlesinger transformations exist for
generic $\{\B_k\}$, see \cite{JM}. Clearly, these transformations change the monodromy data, but they do not change the monodromy matrices $\{M_k\}$ and the Stokes multipliers $S^\pm$.

Now let us take a difference equation of the type considered earlier:
$$
Y(z+1)=A(z)Y(z),\qquad A(z)=A_0z^n+A_1z^{n-1}\dots+A_n.
\tag 5.5
$$
We distinguish two cases (cf. Propositions 1.1 and 1.2): 

\noindent
$\bullet$ $A_0$ is diagonal and has pairwise distinct nonzero eigenvalues;

\noindent
$\bullet$ $A_0=I$, $A_1$ is diagonal and no two eigenvalues of $A_1$ are different by an integer.

As was explained in \S4, see Proposition 4.1, we can generically represent $A(z)$ in the form
$$
A(z)=A_0(z-C_1)\cdots (z-C_n),
$$
where the eigenvalues $\bigl\{a_j^{(i)}\bigr\}$ of $\{C_i\}$ are zeros of $\det A(z)$ divided into $n$ groups of $m$ numbers. We assume, as usual, that no two eigenvalues are different by an integer.

Suppose that $A(z)$ depends on a small parameter $\epsilon$, and as $\epsilon\to 0$, we have
$$
C_i-y_i\epsilon^{-1}+\B_i\to 0,\quad i=1,\dots,n,\qquad
(A_0-I)\epsilon^{-1}-\B_{\infty}\to 0,\qquad \epsilon\to 0,
\tag 5.6
$$
for some pairwise distinct complex numbers $y_1,\dots,y_n$ and some 
$\B_1,\dots,\B_n,\B_\infty\in\Mat$. (The limit relation for $A_0$ is omitted in the case $A_0=I$.)

Note that if we  multiply the unknown function $Y(z)$ in \tht{5.5} by
$\prod_i\Gamma(z-y_i\epsilon^{-1})$, then \tht{5.5} takes the form
$$
\gathered
\left(\prod_i\Gamma(z+1-y_i\epsilon^{-1})Y(z+1)\right)
=(I+B_\infty\epsilon+o(\epsilon))\\ \times
\left(I+\frac{\B_1+o(1)}{z-y_1\epsilon^{-1}}\right)\cdots
\left(I+\frac{\B_n+o(1)}{z-y_n\epsilon^{-1}}\right)
\left(\prod_i\Gamma(z-y_i\epsilon^{-1})Y(z)\right).
\endgathered
$$
If we now assume that $\prod_i\Gamma((\ze+y_i)\epsilon^{-1})Y(\ze\epsilon^{-1})$ tends to a holomorphic function $\Y(\ze)$, then the difference equation above in the limit $\epsilon\to 0$ turns into the differential equation \tht{5.1} (or \tht{5.2}) with $x_i=y_i$.

Substituting the asymptotic relations \tht{5.6} into \tht{4.2}, we see that for any fixed $l_1,\dots,l_n\in\Z$, we have
$$
C_i(l_1,\dots,l_n)+l_i-y_i\epsilon^{-1}+\B_i\to 0,\quad i=1,\dots,n,\qquad \epsilon\to 0.
$$
(This conclusion is based on the fact that if $X=x\epsilon^{-1}+X_0+o(1)$, $Y=y\epsilon^{-1}+Y_0+o(1)$, where $x,y\in\C$, $x\ne y$, and $(z-X)(z-Y)=(z-S)(z-T)$ with $Sp(S)=Sp(Y)$, $Sp(T)=Sp(X)$, then
$S=Y+o(1)$, $T=X+o(1)$, see the explicit construction of $T$ in Lemma 3.3.)

In particular, for any $k_1,\dots,k_n\in\Z$
$$
B_i(k_1,\dots,k_n)+k_i-y_i\epsilon^{-1}+\B_i\to 0,\quad i=1,\dots,n,\qquad \epsilon\to 0.
$$
(See \S4 for the relation of $\{B_i\}$ and $\{C_i\}$.) Thus, on finite intervals $B_i(k)+k_i-y_i\epsilon^{-1}$ is approximately constant.
However, if we assume that $\{B_i(k)+k_i-y_i\epsilon^{-1}\}$ for $k$ of size $\epsilon^{-1}$ approach some smooth functions of $\epsilon k_j$:
$$
\gathered
B_i\bigl([x_1\epsilon^{-1}],\dots,[x_n\epsilon^{-1}]\bigr)+(x_i-y_i)\epsilon^{-1}+\B_i(x_1,\dots,x_n)\to 0,\qquad \epsilon\to 0,\\
\B_i(0,\dots,0)=\B_i,\qquad i=1,\dots,n,
\endgathered
$$
then the corresponding equation \tht{5.5} converges  to \tht{5.1} with $\{\B_i=\B_i(x)\}$ and $x_i$ replaced by $y_i-x_i$. Furthermore, the difference Schlesinger equations \tht{3.1}-\tht{3.3} tend to
$$
\frac{\partial \B_l}{\partial x_j}=\frac
{[\B_l,\B_j]}{(y_j-x_j)-(y_l-x_l)}\,,\qquad   \sum_{l=1}^n
\frac{\partial \B_l}{\partial
x_j}=[\B_l,\B_\infty]\,, \qquad j,l=1,\dots,n.
$$
Comparing these equations to \tht{5.4}, we are led to the following
\proclaim{Conjecture 5.1} For generic $\B_1,\dots,\B_n,\B_\infty$ and pairwise distinct $y_1,\dots,y_n\in\C$, take $B_i=B_i(\epsilon)\in\Mat$, $i=1,\dots,n$ such that
$$
B_i(\epsilon)-y_i\epsilon^{-1}+\B_i\to 0,\qquad \epsilon\to 0.
$$
Let $B_i(k_1,\dots,k_n)$ be the solution of the difference Schlesinger equations \tht{3.1}-\tht{3.3} with the initial conditions $\{B_i(0)=B_i\}$, and let $\B_i(x_1,\dots,x_n)$ be the solution of the classical Schlesinger equations \tht{5.4} with the initial conditions
$$
\B_i(y_1,\dots,y_n)=\B_i, \qquad i=1,\dots,n.
$$
 Then for any $x_1,\dots,x_n\in\Bbb R$ and $i=1,\dots,n$, we have
$$
B_i\bigl([x_1\epsilon^{-1}],\dots,[x_n\epsilon^{-1}]\bigr)+
(x_i-y_i)\epsilon^{-1}+\B_i(y_1-x_1,\dots,y_n-x_n)\to 0,\qquad \epsilon\to 0.
$$
\endproclaim

As for isomonodromy deformations of the second kind (Schlesinger transformations), we are able to prove an asymptotic result rigorously. We will consider the case of equation \tht{5.1}, for \tht{5.2} the situation is similar.

Fix 
$$
\bigl\{t_j^{(k)}\bigr\}_{1\le k\le n,\,\, 1\le j\le m}\subset \C,\qquad 
t_{j_1}^{(k_1)}- t_{j_2}^{(k_2)}\notin\Z\quad\text{unless}\quad j_1=j_2,\, k_1=k_2.
$$
For any 
$$
\B_k\subset\Mat, \quad Sp(\B_k)=\bigl\{t_j^{(k)}\bigr\}_{j=1}^m,\qquad k=1,\dots,n,
$$
and pairwise distinct $x_1,\dots,x_n\in\C$, we define
$$
B_k(\epsilon)=x_k\epsilon^{-1}-\B_k,\qquad k=1,\dots,n,\quad \epsilon\ne 0.
$$
Then we have
$$
Sp(B_k(\epsilon))=\left\{a_j^{(k)}\right\}_{j=1}^m, \quad a_j^{(k)}=x_k\epsilon^{-1}-t_j^{(k)},\qquad j=1,\dots,m.
$$
We also fix $\B_\infty=\diag(s_1,\dots,s_m)$, $s_i\ne s_j$ for $i\ne j$.
Set 
$$
A_0(\epsilon)=I+\epsilon B_\infty.
$$
\proclaim{Lemma 5.2} For generic $\B_1,\dots,\B_n$ and $|\epsilon|$ small enough, there exists a unique degree $n$ polynomial $A(z,\epsilon)=A_0(\epsilon)z^n+A_1(\epsilon)z^{n-1}+\dots$ having $\{z-B_k(\epsilon)\}$ as its right divisors.
\endproclaim
\demo{Proof}
According to Lemma 3.4, the statement is true for large $|\epsilon|$. On the other hand, for fixed $\{\B_k\}$ the existence of $A(z)$ is an open condition on $\epsilon$, and if it holds for large $|\epsilon|$, it also holds for $|\epsilon|$ small enough.\qed
\enddemo

\proclaim{Theorem 5.3} Fix any integers $\bigl\{\lambda_j^{(i)}\bigr\}_{j=1}^m$, $i=1,\dots,n,\infty$, of total sum 0:
$$
\sum_{j=1}^m\left(\sum_{i=1}^n\lambda_j^{(i)}+\lambda_j^{(\infty)}\right)=0.
$$
Then for generic $\B_1,\dots\B_n$ and small enough $|\epsilon|$, there exists the transformation of Theorem 2.1 for the equation $Y(z+1)=A(z,\epsilon)Y(z)$ with 
$$
\kappa_j^{(i)}=-\lambda_j^{(i)},\quad 1\le k\le n,\qquad
\delta_j=-\lambda_j^{(\infty)},\qquad 1\le j\le m.
$$ 
Furthermore, if we denote by $\{\wt \B_k\}$ the coefficients of \tht{5.1} after the corresponding Schlesinger transformation, 
and by $\{\wt B_k(\epsilon)\}$ the matrices such that
$\{z-\wt B_k(\epsilon)\}$ are the right divisors of the transformed $\wt A(z,\epsilon)$, then
$$
\Delta\,\wt B_i(\epsilon)\Delta^{-1}-x_i\epsilon^{-1}+\wt \B_i\to 0,\quad \epsilon\to 0,\quad i=1,\dots,n,
$$
where $\Delta=\diag(\epsilon^{\delta_1},\dots,\epsilon^{\delta_n})$.
\endproclaim
\example{Remark 5.4}
It is easy to verify the statement of Theorem 5.3 if $\lambda_j^{(i)}=-\lambda_j^{(\infty)}=\pm 1$ for some fixed $i$ and all $j=1,\dots,m$, with all other $\lambda$'s being zero. Then 
$\wt B_i(\epsilon)=B_i(\epsilon;0,\dots,0,{ \pm\overset{(i)}\to 1},0,\dots,0).$
As was mentioned above, for any fixed $k_1,\dots,k_n\in\Z$, we have the asymptotics
$$
B_i(\epsilon;k_1,\dots,k_n)+k_i-x_i\epsilon^{-1}+\Cal B_i\to 0,\qquad \epsilon\to 0.
$$
Hence, Theorem 5.3 implies that 
$$
\wt B_k=\cases B_k,& k\ne i,\\
               B_i\pm I, &k=i.
	\endcases
$$	          
This is immediately verified using the fact that the multiplier $\R(\ze)$ for the (continuous) Schlesinger transformation in this case equals $\R(\ze)=(\ze-x_i)^{\pm 1}$.
\endexample
\demo{Proof of Theorem 5.3} Arguing as in the proof of Lemma 5.2, we can
show that all the statements used in this proof which hold generically (like the existence of the polynomial with given right divisors), also hold for generic $\B_1,\dots,\B_n$ and small enough $\epsilon$. Thus, we will ignore the questions of genericity from now on.

Note that we can
decompose the transformations of Theorem 5.3 in both discrete and continuous cases into compositions of elementary ones of the same type (those, for which exactly one of $\{\lambda_j^{(i)}\}$ is equal to $\pm 1$, and exactly one of $\{\lambda_j^{(\infty)}\}$ is equal to $\mp 1$, with all others being zero). It is clear that the claim of the theorem follows from a slightly more general claim for the elementary transformations: we assume that
$$
\Delta_0\,B_l(\epsilon)\Delta_0^{-1}-x_l\epsilon^{-1}+\B_l\to 0,\qquad \epsilon\to 0,\quad l=1,\dots,n,
\tag 5.7
$$
with some diagonal $\Delta_0$ containing integral powers of $\epsilon$ on the diagonal, and we need to conclude that
$$
\Delta\Delta_0\,\wt B_l(\epsilon)\Delta_0^{-1}\Delta^{-1}-x_l\epsilon^{-1}+\wt\B_l\to 0,\qquad \epsilon\to 0,\quad l=1,\dots,n.
\tag 5.8
$$

Let us consider the elementary transformations with $\lambda^{(1)}_j=1$, $\lambda^{(\infty)}_i=-1$.
Denote by $R(z,\epsilon)$ and $\R(\ze)$ the corresponding multipliers for the discrete and continuous equations.
According to the proof of Theorem 2.1, 
$$
R(z,\epsilon)=(z-x_1\epsilon^{-1}+t_j^{(1)})E_i+R_0(\epsilon),\qquad
R^{-1}(z,\epsilon)=I-E_i+\frac{R_1(\epsilon)}{z-x_1\epsilon^{-1}
+t_j^{(1)}}
$$ 
are given by the formulas of Lemma 2.4 with $Q=\hat Y_1(\epsilon)$, and $v=v(\epsilon)$ being an
eigenvector of $B_1(\epsilon)$ with the eigenvalue $x_1\epsilon^{-1}-t_j^{(1)}$. Similarly, \cite{JM, Appendix A} shows that
$$
\R(\ze)=(\ze-x_1)E_i+\R_0,\qquad
\R^{-1}(z)=I-E_i+\frac{\R_1}{\ze-x_1}
$$ 
are given by the same formulas with $Q=\hat \Y_1$ and $v$ being an eigenvector of $\B_1$ with the eigenvalue $t^{(1)}_j$. (Note that only the off-diagonal elements of $Q$ participate in the formulas.) 

\proclaim{Lemma 5.5} Under the assumption \tht{5.7}, we have
$$
\Delta\Delta_0R_0(\epsilon)\Delta_0^{-1}\to \R_0,\qquad \epsilon\,\Delta_0R_1(\epsilon)\Delta_0^{-1}\Delta^{-1}\to \R_1, \qquad\epsilon\to 0
$$
where $\Delta=\epsilon^{E_i}$.
\endproclaim
\demo{Proof} First we note that \tht{5.7} implies that in the projective space the vector $\Delta_0v(\epsilon)$ tends to $v$ as $\epsilon\to 0$. Next, it is easy to see that the difference Schlesinger equations preserve the asymptotics \tht{5.7}: 
$$
\Delta_0\,B_l(\epsilon;k_1,\dots,k_n)\Delta_0^{-1}+k_l-x_l\epsilon^{-1}+\B_l\to 0,\qquad \epsilon\to 0,
\tag 5.9
$$
for any $k_1,\dots,k_n\in\Z$. In particular,
$$
\Delta_0\,C_l(\epsilon)\Delta_0^{-1}-x_l\epsilon^{-1}+\B_l\to 0,\qquad \epsilon\to 0.
$$
Thus, from
$$
A(z,\epsilon)=A_0(\epsilon)(z-C_1(\epsilon))\cdots(z-C_n(\epsilon))=A_0(\epsilon) z^n+A_1(\epsilon)z^{n-1}+\dots
$$
we conclude that
$$
A_1(\epsilon)=-\sum_{l=1}^nC_l(\epsilon)=\text{diagonal part}+\Delta_0^{-1}\left(\sum_{l=1}^n\B_l+o(1)\right)\Delta_0.
$$
By \tht{1.3} we also know that $\epsilon(s_j-s_i){(\hat Y_1(\epsilon))}_{kl}={(A_1(\epsilon))}_{kl}$ for all $k\ne l$. Since $\hat \Y_1=\sum_{l=1}^n\B_i$, the statement follows from the explicit formulas of Lemma 2.4.\qed
\enddemo

A direct computation shows that (here we use the fact that $\R_0\B_1\R_1=t_j^{(1)}\R_0\R_1=0$, which follows from the explicit formulas of Lemma 2.4)
$$
\gather
\wt\B_1=\R_0\left(\B_{\infty}+\sum_{k=2}^{n}\frac{\B_k}{x_1-x_k}\right)\R_1
+E_i\R_1,
\tag 5.10
\\
\wt \B_l=((x_l-x_1)E_i+\R_0)\B_l\left(I-E_i+\frac{\R_1}{x_l-x_1}\right),\quad l=2,\dots,n.
\tag 5.11
\endgather
$$

Let us prove \tht{5.8} for $l\ge 2$ first. Consider the composition of the elementary transformation (for the difference equation) in question with $\s(0,\dots,0,\overset{(l)}\to 1,0,\dots,0)$, see \S3 for the notation. By the uniqueness part of Theorem 2.1, this is equivalent to making $\s(0,\dots,0,\overset{(l)}\to 1,0,\dots,0)$ first, and applying the elementary transformation of after that. Denote the multiplier of this second elementary transformation by $\hat R(z,\epsilon)$. We have
$$
\hat R(z,\epsilon)=(z-x_1\epsilon^{-1}+t_j^{(1)})E_i+\hat R_0(\epsilon),\qquad
\hat R^{-1}(z,\epsilon)=I-E_i+\frac{\hat R_1(\epsilon)}{z-x_1\epsilon^{-1}
+t_j^{(1)}}\,.
$$
 Using Proposition 3.6, we obtain
$$
(z-\wt B_l(\epsilon))R(z,\epsilon)=\hat R(z,\epsilon)(z-B_l(\epsilon)).
$$
Substituting $z=x_l\epsilon^{-1}$ and conjugating by $\Delta_0$, we get
$$
\multline
\Delta_0\wt B_l(\epsilon)\Delta_0^{-1}-x_l\epsilon^{-1}\\=
\left(\Delta_0 \hat R(x_l\epsilon^{-1},\epsilon)\Delta_0^{-1}\right)
\left(\Delta_0 B_l(\epsilon)\Delta_0^{-1}-x_l\epsilon^{-1}\right)\,
\left(\Delta_0  R^{-1}(x_l\epsilon^{-1},\epsilon)\Delta_0^{-1}\right).
\endmultline
$$
Because of \tht{5.9}, the limit relations of Lemma 5.5 also hold for $\hat R_0$, $\hat R_1$. Using them, Lemma 5.5 itself, \tht{5.7} and \tht{5.11}, we arrive at \tht{5.8} for $l\ge 2$.

Thus, it remains to prove \tht{5.8} for $l=1$. We have
$$
\gathered
\wt A(z,\epsilon)=R(z+1,\epsilon)A_0(z-C_1(\epsilon))\cdots(z-C_n(\epsilon))
R^{-1}(z,\epsilon)\\=A_0(z-\wt C_1(\epsilon))\cdots(z-\wt C_n(\epsilon)).
\endgathered
\tag 5.12
$$

The relation \tht{5.8} for $l\ge 2$ implies that 
$$
\Delta\Delta_0\,\wt C_l(\epsilon)\Delta_0^{-1}\Delta^{-1}-x_l\epsilon^{-1}+\wt\B_l\to 0,\qquad \epsilon\to 0,\quad l=2,\dots,n.
$$

Substituting these estimates and similar ones for $C_l$ and setting $z=(w+x_1)\epsilon^{-1}-t_j^{(1)}$, we can rewrite \tht{5.12} as follows (note that $A_0$ is diagonal and hence it commutes with $\Delta_0$):
$$
\gathered
((w-x_1+\epsilon)E_i+\Delta\Delta_0R_0(\epsilon)\Delta_0^{-1})
(I+\epsilon\B_\infty)\left(I+\epsilon\,\frac{\B_1+o(1)}{w-x_1}\right)
\cdots \left(I+\epsilon\,\frac{\B_{n}+o(1)}{w-x_{n}}\right)\\
\times \left(I-E_i+\frac{\epsilon\Delta_0R_1(\epsilon)\Delta_0^{-1}\Delta^{-1}}
{w-x_1}\right)=(I+\epsilon\B_\infty)
\left(I+\epsilon\,\frac{x_1\epsilon^{-1}-\Delta\Delta_0\wt C_1(\epsilon)\Delta_0^{-1}\Delta^{-1}}{w-x_1}\right)\\ \times
\left(I+\epsilon\,\frac{\wt\B_2+o(1)}{w-x_2}\right)
\cdots \left(I+\epsilon\,\frac{\wt\B_{n}+o(1)}{w-x_{n}}\right).
\endgathered
$$
Comparing the residues of both sides at $w=x_1$ and looking at terms of order $\epsilon$, we see that
$$
x_1\epsilon^{-1}-\Delta\Delta_0\wt C_1(\epsilon)\Delta_0^{-1}\Delta^{-1}\to\wt B_1,
$$
where $\wt B_1$ is given by \tht{5.10}.
(We need to use Lemma 5.5 and the relation $R_0(\epsilon)R_1(\epsilon)=0$ here.) Since the difference Schlesinger equations preserve the asymptotics \tht{5.8} (cf. \tht{5.9}), we get \tht{5.8} for $l=1$, and thus for all $l$.

The proof that \tht{5.7} implies \tht{5.8} in the case $\lambda_j^{(1)}=-\lambda_j^{(\infty)}=-1$ is very similar. Let us outline the necessary changes. The multipliers have the form
$$
R(z,\epsilon)=I-E_i+\frac{R_1^t(\epsilon)}{z-1-x_1\epsilon^{-1}
+t_j^{(1)}},\qquad
R^{-1}(z,\epsilon)=
(z-1-x_1\epsilon^{-1}+t_j^{(1)})E_i+R_0^t(\epsilon),
$$ 
where $R_0(\epsilon)$, $R_1(\epsilon)$ are constructed by Lemma 2.4 with 
$v=v(\epsilon)$ being a solution of $A^t(x_1\epsilon^{-1}-t_j^{(1)},\epsilon)\,v(\epsilon)=0$ and $Q=-\hat Y_1^t(\epsilon)$, see the proof of Theorem 2.1. Similarly,
$$
\R(\ze)=I-E_i+\frac{\R_1^t}{\ze-x_1},\qquad
\R^{-1}(\ze)=
(\ze-x_1)E_i+R_0^t,
$$
where $\R_0$ and $\R_1$ are as in Lemma 2.4 with $v$ being an eigenvector
of $\B_1^t$ with the eigenvalue $t_j^{(1)}$, and $Q=-\hat \Y_1^t$.
Similarly to Lemma 5.5, \tht{5.7} implies
$$
\Delta_0R_0^t(\epsilon)\Delta_0^{-1}\Delta^{-1}\to \R_0^t,\qquad \epsilon\,\Delta\Delta_0R_1^t(\epsilon)\Delta_0^{-1}\to \R_1^t, \qquad\epsilon\to 0
$$
with $\Delta=\epsilon^{-E_i}$. Similarly to \tht{5.10}, \tht{5.11}, we have
$$
\gather
\wt\B_1=\R_1^t\left(\B_{\infty}+\sum_{k=2}^{n}\frac{\B_k}{x_1-x_k}\right)\R_0^t
-\R_1^tE_i,
\\
\wt \B_l=\left(I-E_i+\frac{\R_1^t}{x_l-x_1}\right)\B_l\,((x_l-x_1)E_i+\R_0^t),\quad l=2,\dots,n.
\endgather
$$

Using the same argument with composing our elementary transformation with $\s(0,\dots,0,\overset{(l)}\to 1,0,\dots,0)$, we prove \tht{5.8} for $l\ge 2$. Then substituting estimates for $C_l$'s and $\wt C_l$'s into 
$\wt A(z,\epsilon)=R(z+1,\epsilon)A(z,\epsilon)R^{-1}(z,\epsilon)$, we get (with $z=(w+x_1)\epsilon^{-1}-t_j^{(1)}$)
$$
\gathered
\left(I-E_i+\frac{\epsilon\Delta\Delta_0R_1^t(\epsilon)\Delta_0^{-1}}
{w-x_1}\right)(I+\epsilon\B_\infty)\left(I+\epsilon\,\frac{\B_1+o(1)}{w-x_1}\right)
\cdots \left(I+\epsilon\,\frac{\B_{n}+o(1)}{w-x_{n}}\right)\\
\times ((w-x_1-\epsilon)E_i+\Delta_0R_0^t(\epsilon)\Delta_0^{-1}\Delta^{-1})=(I+\epsilon\B_\infty)\\ \times
\left(I+\epsilon\,\frac{x_1\epsilon^{-1}-\Delta\Delta_0\wt C_1(\epsilon)\Delta_0^{-1}\Delta^{-1}}{w-x_1}\right)
\left(I+\epsilon\,\frac{\wt\B_2+o(1)}{w-x_2}\right)
\cdots \left(I+\epsilon\,\frac{\wt\B_{n}+o(1)}{w-x_{n}}\right).
\endgathered
$$
Comparing the residues of both sides at $w=x_1$ and taking terms of order $\epsilon$, we recover the estimate of type \tht{5.8} for $\wt C_1(\epsilon)$, and hence for $\wt B_1(\epsilon)$. The proof of Theorem 5.3 is complete.\qed 
\enddemo

\end